\documentclass[a4,11pt]{article}

\usepackage[latin1]{inputenc}  %% les accents dans le fichier.tex
\usepackage[T1]{fontenc}       %% Pour la c\'esure des mots accentu\'es
\usepackage{vmargin}
\usepackage{amsmath,amssymb}
\usepackage{epsfig,color}
\usepackage[round]{natbib}

\newcommand{\1}{1 \hspace*{-0.2ex}\rule{0.10ex}{1.5ex}\hspace{0.2ex}}

\newcommand{\E}{\mathbb{E}}
\newcommand{\N}{\mathbb{N}}

\renewcommand{\P}{\mathbb{P}}
\newcommand{\R}{\mathbb{R}}
\newcommand{\Z}{\mathbb{Z}}
\def\mA{\mathcal A}
\def\tA{\tilde{\mathcal A}}
\def\mE{\mathcal E}
\def\mR{\mathcal R}
\def\mB{\mathcal B}
\def\ds{\displaystyle}

\newtheorem{Theorem}{Theorem}[section]
\newtheorem{Theoreme}[Theorem]{Theorem}

\newtheorem{Proposition}[Theorem]{Proposition}

\newtheorem{Remark}[Theorem]{Remark}

\newtheorem{Hypothesis}[Theorem]{Hypothesis}

\def\un{{\mathrm{1~\hspace{-1.4ex}l}}}
\def\ds{\displaystyle}

\title{Convergence of stochastic gene networks to hybrid piecewise deterministic processes}

\author{ A. Crudu$^{1}$, A. Debussche$^{2}$, A. Muller$^{3}$, O. Radulescu$^{4}$ \\
\small
$^{1}$ IRMAR - UMR 6625, Universit\'e de Rennes 1, Campus de Beaulieu, \\
\small 35042 Rennes, France, \\
\small
$^{2}$ IRMAR - UMR 6625, ENS Cachan Bretagne, Campus de Ker Lann, \\
\small 35170 Bruz, France, \\
\small$^{3}$ IECN - UMR CNRS 7502, Universit\'e Henri Poincar\'e\\
\small 54506  Vandoeuvre-l\`es-Nancy, France\\
\small
$^{4}$DIMNP - UMR 5235 CNRS/UM1/UM2, Universit\'e de Montpellier 2, Place Eug\`ene Bataillon, \\ \small CP 107,
34095 Montpellier, France.
}

\begin{document}
\maketitle

\abstract{We study the asymptotic behavior of multiscale stochastic gene networks using
weak limits of Markov jump processes. Depending on the time and concentration scales of
the system we distinguish four types of limits: continuous piecewise deterministic processes (PDP)
with switching, PDP with jumps in the continuous variables, averaged PDP, and PDP with singular switching.
We justify rigorously the convergence for the four types of limits. The convergence results can be used
to simplify the stochastic dynamics of gene network models arising in molecular biology.
}


\section{Introduction}

Modern molecular biology emphasizes the important role of the gene regulatory networks in the
functioning of living organisms. Recent experimental advances in molecular biology show
that many gene products do not follow deterministic dynamics
and should be modeled as random variables (\cite{kepler2001stochasticity,kaufmann2007stochastic}).

Markov processes approaches to gene networks dynamics, originating from the pioneering ideas
of Delbr\"uck (\cite{delbruck}), capture diverse features of the
experimentally observed expression variability,
such as bursting (\cite{cai06}), various types of steady-state
distributions of RNA and protein numbers (\cite{kaern}),
noise amplification or reduction by network propagation (\cite{paulsson,warren2006ern}),
clock de-synchronization (\cite{barkai2000biological}),
stochastic transitions in cellular memory storage circuits (\cite{oudernaden07}).

However, the study of the full Markov dynamics
of biochemical networks is a difficult task. Even the simplest Markovian model, such as the
production module
of a single protein involves tens of variables and biochemical reactions and an equivalent
number of parameters (\cite{kierzek01,krishna}). The direct simulation of such models by the
Stochastic Simulation Algorithm (SSA) (\cite{gillespie76}) is extremely time consuming.
%Stochastic
%simulation in molecular biology is a field of modeling that strongly needs simplified models.

In order to increase computational efficiency, several accelerated
simulation algorithms are hybrid and treat fast biochemical reactions as continuous
variables (\cite{haseltine,alfonsi2005ash,alfonsi}). Similar approaches reducing fast
reactions can be justified by diffusion approximations for Markov processes (\cite{ball06}).

A different hybrid approach is to distinguish between molecular species according to their abundances.
Species in small amounts can be treated as discrete variables, whereas species in large amounts
can be considered continuous. It has been proposed that,
the dynamics of gene networks with well separated abundances,
can be well approximated by  piecewise deterministic Markov processes (\cite{TSI,crudu2009hybrid}).
Piecewise deterministic processes (PDP) are used in operational research in relation with
optimal control and various technological applications (\cite{boxma,ghosh_msh,pola2003shm,bujorianu2004gsh}).
Their popularity
in physical, chemical  and biological sciences is also steadily increasing as
they provide a natural framework to deal with intermittent phenomena in many contexts (\cite{TSI,zeiser2008sgn}).

By looking for the best PDP approximation of a stochastic network of biochemical reactions,
and depending on the time scales of the reaction mechanism, we can distinguish several
cases (\cite{crudu2009hybrid}):

\begin{itemize}
\item Continuous PDP with switching:  continuous variables evolve according to ordinary differential equations. The trajectories of
the continuous variables are continuous, but the differential equations depend on one or several discrete variables.
\item PDP with jumps in the continuous variables: the same as the previous case, but the continuous variables can jump
as well as the discrete variables.
\item Averaged PDP: some discrete variables have rapid transitions and can be averaged. The resulting approximation is an averaged PDP.
\item Discontinuous PDP with singular switching: the continuous variable has two time scaling. The switch between the two regimes
is commanded by a discrete variable. The rapid parts of the trajectory of the continuous variable can be approximated by discontinuities.
\end{itemize}

In this paper, we justify rigorously these approximations that
were illustrated by models of stochastic gene expression in \cite{crudu2009hybrid}.
More precisely, we present several
theorems on the weak convergence of biochemical reactions processes  towards
piecewise deterministic processes of the type specified above.
The resulting piecewise deterministic processes can be used for more efficient simulation
 algorithms, also, in certain
cases, can lead to analytic results for the stochastic behavior of gene networks.
Higher order approximations of multiscale stochastic chemical kinetics, corresponding
to stochastic differential equations with jumps, though not discussed in this paper,
represent straightforward extensions of our results.

The structure of this article is as follows. In section 2, we present the PDP, a useful theorem on the uniqueness of
the solution of a martingale problem and the Markov jump model for stochastic regulatory networks.
The four remaining sections discuss the asymptotic
behaviors of the models, corresponding to the four cases presented above.

\section{Piecewise Deterministic Processes}

We begin with a brief description of  Piecewise Deterministic Processes (PDP) and collect useful results on these. We do not consider PDPs in their full generaliy. The reader
is refered to \cite{Davis} for further results.

\bigskip
\noindent{\textbf{Standard conditions: }}

In this article, a PDP taking values in $E=\R^n\times \N^d$ is a process $x_t=(y_t, \nu_t)$, determined by its three local characteristics :
\begin{enumerate}
\item For all $\nu\in \N^d$, a Lipschitz continuous vector field in $\R^n$, denoted by $F_\nu$, {which} determines a unique global flow $\phi_\nu(t,y)$ in $\R^n$ such that, for $t>0$,
$$\frac{d}{dt}\phi_\nu(t,y)=F_\nu(\phi_\nu(t,y)),\qquad \phi_\nu(0,y)=y,\qquad \forall y\in \R^n.$$ We also use the notation: $F(y,\nu)=F_\nu(y)$.

\item A jump rate $\lambda : E\rightarrow \R^+$ such that, for each $x=(y,\nu)\in E$, there exists $\epsilon(x)>0$ such that
$$\int_0^{\epsilon(x)}\lambda(\phi_\nu(t,y),\nu)dt<\infty.$$
\item A transition measure $Q:E\rightarrow\mathcal P(E)$, $x\mapsto Q(\cdot;x)$, where $\mathcal P(E)$ denotes the set of probability measures on $E$. We assume that $Q(\{x\};x)=0$ for each $x\in E$.
\end{enumerate}
From these standard conditions, a right-continuous sample path $\{x_t : t>0\}$ starting at $x=(y,\nu)\in E$ may be constructed as follows. Define
$$x_t(\omega):=\phi_{\nu}(t,y),\qquad \textrm{ for } 0\leq t< T_1(\omega),$$
where $T_1(\omega)$ is the realization of the first jump time $T_1$, with the following distribution~:
$$\P_x(T_1>t)=\exp\Big(-\int_0^t\lambda(\phi_\nu(s,y),\nu)ds\Big)=:H(t,x),\qquad t\in \R^+.$$
We have then $x_{_{T_1^-(\omega)}}(\omega)=(\phi_\nu(T_1(\omega),y),\nu)$, and the post-jump state $x_{T_1(\omega)}(\omega)$ has the distribution given by :
$$\P_x(x_{T_1}\in A| T_1=t)=Q(A;(\phi_\nu(t,y),\nu))$$
on the Borel sets $A$ of $E$.

We then restart the process at $x_{T_1(\omega)}(\omega)$ and proceed recursively according to the same procedure to obtain a sequence of jump-time realizations $T_1(\omega),T_2(\omega),\ldots$.
 Between each of two consecutive jumps, $x_t(\omega)$ follows a deterministic motion, given by the flow corresponding to the vector field $F$.

Such a process $x_t$ is called a PDP. The number of jumps that occurr between the times 0 and $t$ is denoted by
$$N_t(\omega)=\sum_k\1_{t\geq T_k}(\omega).$$

It can be shown that $x_t$ is a strong Markov process with right-continuous, left-limited sample paths (see \cite{Davis}). The generator $\mathcal A$ of the process is formally given by
\begin{equation}
\label{APDP}
\mathcal Af(x)=F_\nu(x)\cdot \nabla_y f(x)+\lambda(x)\int_{E}(f(z)-f(x))Q(dz;x)
\end{equation}
for each $x=(y,\nu)\in E$, we have denoted by $\nabla_y$ the gradient
with respect to the variable $y\in\R^n$. The  domain of $\mathcal A$
is described precisely in  \cite{Davis}. We do not
need such a precise description and just note that $\mathcal A$ is well defined for
$f\in \mathcal E$,  the set of functions $f:E\to \R$ such that :
\begin{itemize}
\item[E1.] $f$ is bounded,
\item[E2.]for all $\nu\in\N^d$, $f(\cdot,\nu)\in C^1(\R^n)$,
\item[E3.] its derivatives are bounded uniformly in $E$.
 \end{itemize}
For $f\in\mathcal E$, we denote by
\begin{equation}
L_f=\sup_{y\in\R^n}||D_yf||_\infty=\sup_{(y,\nu)\in E}\|D_y f(y,\nu)\|,
\end{equation}
the Lipschitz constant of $f$ with respect to the variable $y$.

The space $\mathcal E$ is a Banach space when endowed with the norm
\begin{equation}
||f||_{\mathcal E}=||f||_{\infty}+L_f.
\end{equation}

\par
\medskip
If $Z$ is a Banach space, $\mathcal B_b(Z)$ is the set of bounded Borel measurable functions on $Z$~; $C_b^k(Z)$ is the set of $C^k$-differentiable functions on $Z$, such that the derivatives, until the $k$-th order, are bounded~; $\mathcal C_b(Z)$ is the set of bounded continuous functions on $Z$. Also $D(\R^+;Z)$ is the set of process defined on $\R^+$ with right-continuous, left-limited sample paths defined on $\R^+$ and taking values in $Z$ and $C(\R^+;Z)$ is the set of continuous process defined on $\R^+$ and taking values in $Z$.

\par
\medskip
The PDPs considered in this paper will always satisfy the following property~:

\begin{Hypothesis}
\label{h1}
The three local characteristics of the PDP satisfy the standard conditions given above. The jump rate $\lambda$ is $C^1$-differentiable with respect to the variable $y\in \R^n$. For every starting point $x=(y,\nu)\in E$ and $t\in \R^+$, we suppose $\E(N_t)<\infty$.
\end{Hypothesis}

\begin{Remark}
$\E(N_t)<\infty$ implies in particular that $T_k(\omega)\to\infty$ almost surely. This assumption is usually quite easy to check in applications, but it is hard to formulate general conditions under which it holds, because of the complicated interaction between $F,\lambda$ and $Q$. It can be shown for instance that if $\lambda$ is bounded, then $\E(N_t)<\infty$ (cf. \cite{Davis}).
\end{Remark}

For some results, we need the following stronger property~:
\begin{Hypothesis}
\label{h2}
The functions  $F$, $\lambda$ and $x\mapsto \lambda(x)\int_E f(z) Q(dz;x)$,
with $f\in \mE$, are bounded on $E$, $C^1$-differentiable with respect to the variable $y\in \R^n$ and their derivatives with respect to $y$ are also bounded.
\end{Hypothesis}
When Hypothesis \ref{h2} is satisfied, we set
$$
\begin{array}{l}
\ds M_F=\|F\|_\infty, \quad L_F=\sup_{y\in\R^n} \|D_yF\|_\infty,\\
\\
\ds M_\lambda = \|\lambda\|_\infty,\quad L_\lambda = \sup_{y\in\R^n} \|D_y\lambda\|_\infty.\\
\end{array}
$$
and $L_Q$ a constant such that, for all $f\in \mE$ and $x=(y,\nu)\in E$ :
$$
\left\| D_y \left(\lambda(x) \int_E f(z) Q(dz;x)\right)\right\|_\infty\le L_Q\|f\|_\mE
$$
Denote by $(P_t)_{t\ge 0}$ the transition semigroup associated to the PDP constructed above
and by $P_x$ the law of the PDP starting from $x\in E$. Then, $P_x$ is a solution of the
martingale problem associated to $\mathcal A$ in the following sense:
$$
f(x_t)-f(x)-\int_0^t {\mathcal A}f(x_s)ds
$$
is a local martingale for any $f\in \mE$ (see \cite{Davis}). Moreover, if Hypothesis \ref{h2} holds, it is a bounded
martingale. As usual, we have denoted by $(x_t)_{t\ge 0}$ the canonical process on $D(\R^+;E)$.

The following results gives a uniqueness property for this martingale problem. It will enable us to characterize the asymptotic behavior of our stochastic regulatory networks.
\begin{Theoreme}
\label{t2.2}
If Hypothesis \ref{h2} is satisfied, then the law of the PDP determined by $F$, $\lambda$, and $Q$ is the unique solution of the martingale problem associated to the generator $\mathcal A$.
\end{Theoreme}

The proof of the theorem is given in the appendix.

\par
\bigskip
\noindent{\bf{Markov jump model for stochastic regulatory networks: known results}}

 We consider a set of chemical reactions $R_r$, $r\in \mR$; $\mR$ is supposed to be finite. These reactions involve species indexed by a set $S={1,\dots,M}$, the number of molecules of the specie
$i$ is denoted by $n_i$ and $X\in\N^{M}$ is the vector consisting of the $n_i$'s.
Each reaction $R_r$  has a rate
$\lambda_r(X)$ which depends on the state of the system, described by $X$ and corresponds
to a change $X\to X+\gamma_r$, $\gamma_r\in\Z^M$.

Mathematically, this evolution can be described by the following Markov jump process. It is based on a
sequence $(\tau_k)_{k\ge 1}$ of random waiting times with exponential distribution. Setting $T_0=0$,
$T_i= \tau_1+\dots +\tau_i$, $X$ is constant on $[T_{i-1}, T_{i})$ and has a jump at $T_{i}$.
The parameter of $\tau_i$ is given
by $\sum_{r\in {\mathcal R}} \lambda_r(X(T_{i-1}))$:
$$
\P(\tau_i>t)= \exp\Big(-\sum_{r\in {\mathcal R}} \lambda_r(X(T_{i-1}))t\Big).
$$
At time $T_i$, a reaction $r\in {\mathcal R}$ is chosen with probability
$\lambda_r(X(T_{i-1}))/\sum_{r\in {\mathcal R}} \lambda_r(X(T_{i-1}))$ and the state changes according
to $X\to X+\gamma_r$:
$$
X(T_{i})= X(T_{i-1})+\gamma_{r}.
$$
This Markov process has the following generator (see \cite{ethierkurtz}):
$$
A f(X) = \sum_{r\in \mR} \left[ f(X+\gamma_r)-f(X)\right]\lambda_r(X).
$$
We do not need a precise description of the domain of $A$, the above definition holds for instance
for functions in $C_b(\R^{M})$.

In the applications we have in mind, the numbers of molecules have different scales. Some of
the molecules are in small numbers and some are in large numbers.  Accordingly, we split
the set of species into two sets $C$ and $D$ with cardinals $M_C$ and $M_D$.
This induces the decomposition
$X=(X_C,X_D)$, $\gamma_r=(\gamma_r^C,\gamma_r^D)$.
For $i\in D$, $n_i$ is of order $1$ while for
$i\in C$, $n_i$ is proportional to $N$ where $N$ is a large number. For $i\in C$, setting $\tilde n_i = n_i/N$, $\tilde n_i$ is of order $1$. We define
$\ds x_C=\frac1N X_C$ and $x=(x_C,X_D)$.

We also decompose the set of reactions according to the species involved. We set
$\mR=\mR_D\cup \mR_C\cup \mR_{DC}$. A reaction in $\mR_D$ (resp. $\mR_C$)
produces or consumes only species in $D$ (resp. $C$). Also, the rate of a reaction
in $\mR_{D}$ (resp. $\mR_{C}$) depends only on $X_{D}$ (resp. $x_{C}$).
A reaction in  $\mR_{DC}$ has a rate depending on both $x_C$ and $X_D$ and produces or consumes, among others, species from
$C$ or $D$.

The rate of a reaction in $r\in\mR_{C}$ is also large and of order $N$ and we set
$\tilde \lambda_{r}= \frac{\lambda_{r}}{N}$. In general, reactions in $\mR_{D}$ or $\mR_{DC}$ have a rate of order $1$.

Introducing the new scaled variables, the generator has the
form:
$$
\begin{array}{ll}
\tilde \mA f(x_{C},X_{D}) &=\ds  \sum_{r\in \mR_{C}} \left[ f(x_{C}+\frac1N \gamma_r^C,X_{D})-f(x_{C},X_{D})
\right]N\tilde\lambda_r(x_{C})\\
\\
&\ds +\sum_{r\in \mR_{DC}} \left[ f(x_{C}+\frac1N\gamma_r^C,X_{D}+\gamma_{r}^D)-f(x_{C},X_{D})\right]\lambda_r(x_{C},X_{D})\\
\\
&\ds +\sum_{r\in \mR_{D}} \left[ f(x_{C}, X_{D}+\gamma_r^D)-f(x_{C},X_{D})\right]\lambda_r(X_{D}).
\end{array}
$$
Assuming that the scaled rates $\tilde \lambda_r$ are $C^1$ with respect to $x_C$, it is not difficult to see that if $N\to\infty$, the two sets of species decouple. Indeed,
reactions in $\mR_{DC}$ do not happen sufficiently often and they do not change $x_{C}$
in a sufficiently large manner. The limit would simply give a set of differential equations for
the continuous variable $x_{C}$, which evolves without influence of $X_D$. The discrete variable would have its own dynamic made of jumps. These results have been shown by \cite{kurtz71} and \cite{kurtz78}.

In the following sections, we consider more general systems where other types of reactions may happen and which yield different limiting systems

\section{Continuous piecewise deterministic process}
\label{CPDP}

In this section, we assume that some of the reactions in a subset $S_{1}$ of $\mR_{DC}$ are such that
their rate is large and scales with $N$. We again set $\tilde \lambda_{r}= \frac1N\lambda_{r}$ for
$r\in S_{1}$. We assume that these equations do not change $X_{D}$, in other words \begin{equation}
\label{e2.3.1}
\gamma_r^D=0,\; r\in S_1.
\end{equation}
 However, the rate $\lambda_r$ depends on $X_D$. Note that this is possible and even frequent in molecular biology, meaning
 that reactions of the type $S_1$ recover the reactant, like in the reaction
 $A \to A + B$, with $A$ and $B$ discrete and continuous species, respectively. 
 The more complicated case $\gamma_r^D\neq 0$ is treated in section 5.

The scaled generator has now the form
\begin{equation}
\label{e1}
\begin{array}{ll}
\tilde \mA_{N} f(x_{C},X_{D}) &=\ds  \sum_{r\in \mR_{C}} \left[ f(x_{C}+\frac1N \gamma_r^C,X_{D})-f(x_{C},X_{D})
\right]N\tilde\lambda_r(x_{C})\\
\\
&\ds +\sum_{r\in S_{1}} \left[ f(x_{C}+\frac1N\gamma_r^C,X_{D})-f(x_{C},X_{D})\right]N\tilde\lambda_r(x_{C},X_{D})\\
\\
&\ds +\sum_{r\in \mR_{DC}\setminus S_{1}} \left[ f(x_{C}+\frac1N\gamma_r^C,X_{D}+\gamma_{r}^D)-f(x_{C},X_{D})\right]\lambda_r(x_{C},X_{D})\\
\\
&\ds +\sum_{r\in \mR_{D}} \left[ f(x_{C}, X_{D}+\gamma_r^D)-f(x_{C},X_{D})\right]\lambda_r(X_{D}).
\end{array}
\end{equation}
%We have added the subscript $N$ to emphasize the dependance on the volume.
For $f\in C^1_{b}(E)$, we may let $N\to \infty$ and obtain  the limit
generator
$$
\begin{array}{ll}
\ds \mA_{\infty}f(x_{c},X_{D}) &=\ds  \left(\sum_{r\in \mR_{C}} \tilde \lambda_{r}(x_{C})  \gamma_{r}^C+\sum_{r\in S_{1}}  \tilde \lambda_{r}(x_{C},X_{D})  \gamma_{r}^C \right)
\cdot \nabla_{x_{C}}f(x_{C},X_{D}) \\
\\
&\ds  +\sum_{r\in \mR_{DC}\setminus S_{1}} \left[ f(x_{C},X_{D}+\gamma_{r}^D)-f(x_{C},X_{D})\right]\lambda_r(x_{C},X_{D})\\
\\
&\ds +\sum_{r\in \mR_{D}} \left[ f(x_{C}, X_{D}+\gamma_r^D)-f(x_{C},X_{D})\right]\lambda_r(X_{D}).
\end{array}
$$
This formal argument indicates that, as $N\to\infty$, the process converges to a continuous PDP (see (\ref{APDP})). The state is described by a continuous variable $x_{C}$ and a discrete
variable $X_{D}$. The discrete variable is a jump process and is piecewise constant.
The continuous variable evolves according to differential equations depending
on $X_{D}$. It is continuous but the vector field describing its evolution changes when $X_{D}$
jumps.

This is rigorously justified by the following theorem.

\begin{Theoreme}
\label{t2.3.1}
Let $x^N=(x_C^N,X^N_D)$ be a jump Markov process as above, starting at $x^N(0)=(x_C^N(0),X_D^N(0))$. Assume that the jump rates $\tilde \lambda_r$, $r\in \mR_C\cup S_1$ and
$ \lambda_r$, $r\in \mR_{DC}\setminus S_1$ are $C^1$ functions of $x_C\in \R^{M_C}$.
We define $P_{x_0}$ the law of the PDP
starting at $x_0=(x_{C,0},X_{D,0})$ whose jump intensities are:
$$\lambda(x)=\sum_{r\in {\mR_D \cup \mR_{DC} \setminus S_1 }} \lambda_r(x),$$
the transition measure is defined by:
$$
\begin{array}{l}
\int_{E}f(z)Q(dz;x)\\
\ds =\frac{1}{\lambda(x)}\left(\sum_{r\in \mR_{DC}\setminus S_1}f(x_C,X_D+\gamma_r^D)\lambda_r(x_C,X_D)+\sum_{r\in \mR_D}f(x_C,X_D+\gamma_r^D)\lambda_r(X_D)
\right),
\end{array}
$$ 
for $x=(x_C,X_D),$
and the vectors fields are given by:$$F_{X_D}(x_C)=\sum_{r\in\mR_C}\gamma_r^C\tilde\lambda_r(x_C)+\sum_{r\in S_1}\gamma_r^C\tilde\lambda_r(x_C,X_D).$$

%, with jump intensities $\lambda_r, r\in \mR$, and transition measure $Q(x+\gamma_r,x)=\lambda_r(x)/\sum_r\lambda_r(x)$.
Assume that Hypothesis \ref{h1} is satisfied and $x^N(0)$ converges in distribution to $x_0$, then $x^N$ converges in distribution to the PDP whose law is $P_{x_0}$.

\end{Theoreme}

{\bf{Proof:}}
\par
In the following, we work only with scaled variables and simplify the notation by omitting the tildes.
In other words, we use $\lambda_r$ to denote the rate of all reactions.
\par
The proof is divided into three steps. We begin our proof by supposing that the jump rates and their derivatives with respect to $x_C$ are bounded. Hypothesis \ref{h2} is then satisfied. We then prove Theorem \ref{t2.3.1} by a truncation argument.
\par
\medskip
{\it{Step 1: Tightness for bounded reaction rates.}}
\par

We first assume that all rates $\lambda_r$ are bounded as well as their derivatives
with respect to $x_C$.

Let $x^N$ be a Markov jump process whose generator is given by $\tilde \mA_{N}$.

Without loss of generality, we assume that the initial value of the process is deterministic:
$x^N(0) = (x^N_{C}(0), X^N_{D}(0))$ and converges  to  $x_{0}=(x_{C,0},X_{D,0})$
in $\R^{M_{C}}\times \N^{M_{D}}$.

Let $(Y_{r})_{r\in \mR}$ be a sequence of independent standard Poisson processes.
By Proposition 1.7, Part 4, and Theorem 4.1, Part 6, of \cite{ethierkurtz} we know that there exists stochastic processes $(\tilde x^N)_{N\in\N}$ in $D(\R^+; E)$ such that
$$
\tilde x^N(t)= x^N(0)+\sum_{r\in\mR} \gamma_{r} Y_{r}\left( \int_{0}^t \lambda_{r}(\tilde x^N(s)) ds \right), \;
t\ge 0.
$$
Moreover, for each $N$, $x^N$ and $\tilde x^N$ have the same distribution. Since we consider only the
distributions of the processes, we only consider $\tilde x^N$ in the following and use the same notation
for both processes.

Using the decomposition $x^N=(x_C^N,X^N_D)$, we have
$$
\begin{array}{ll}
x^N_{C}(t)& \ds = x^N_{C}(0)
+ \sum_{r\in\mR_{C}} \frac1N \gamma_{r}^CY_{r}\left( N\int_{0}^t \tilde \lambda_{r}(x^N_{C}(s))ds
\right) \\
\\
&\ds+ \sum_{r\in S_{1}} \frac1N \gamma_{r}^CY_{r}\left( N\int_{0}^t  \tilde\lambda_{r}(x^N_{C}(s),X^N_{D}(s))ds
\right) \\
\\
&\ds+ \sum_{r\in \mR_{DC}\setminus S_{1}} \frac1N \gamma_{r}^CY_{r}\left( \int_{0}^t  \lambda_{r}(x^N_{C}(s),X^N_{D}(s))ds
\right)
\end{array}
$$
and
$$
\begin{array}{ll}
X^N_{D}(t)& \ds = X^N_{D}(0)
+ \sum_{r\in \mR_{DC}\setminus S_{1}}  \gamma_{r}^D Y_{r}\left( \int_{0}^t  \lambda_{r}(x^N_{C}(s),X^N_{D}(s))ds
\right) \\
\\
&\ds+ \sum_{r\in \mR_{D}}  \gamma_{r}^DY_{r}\left( \int_{0}^t  \lambda_{r}(X^N_{D}(s))ds.
\right)
\end{array}
$$
We easily prove tightness
in $D(\R^+;\R^{M_D})$ of the laws of $(X^N_D)_{N\in\N}$ by the same proof as for Proposition 3.1 in chapter 6 of \cite{ethierkurtz} and by using the fact that the law of $Y_r$ is tight in $D(\R^+;\N)$, for every $r\in \mR$, according to Theorem 1.4 of \cite{bill}.

To prove that the laws of  $(x^N_C)_{N\in\N}$ are tight in $C(\R^+; \R^{M_C})$, we adapt the
argument of section 2 chapter 11 in \cite{ethierkurtz}.

Let $\tilde Y_r(u) =Y_r(u) - u$ be the standard Poisson process centered at its expectation, we have:
$$
\begin{array}{ll}
x^N_{C}(t)& \ds = x^N_{C}(0)
+ \sum_{r\in\mR_{C}} \frac1N \gamma_{r}^C\tilde Y_{r}\left( N\int_{0}^t  \tilde\lambda_{r}(x^N_{C}(s))ds
\right) \\
\\
&\ds+ \sum_{r\in S_{1}} \frac1N \gamma_{r}^C\tilde Y_{r}\left( N\int_{0}^t \tilde\lambda_{r}(x^N_{C}(s),X^N_{D}(s))ds
\right) \\
\\
&\ds+\int_0^t F(x^N_C(s),X^N_D(s))ds\\
\\
&\ds+ \sum_{r\in \mR_{DC}\setminus S_{1}} \frac1N \gamma_{r}^CY_{r}\left( \int_{0}^t  \lambda_{r}(x^N_{C}(s),X^N_{D}(s))ds
\right)
\end{array}
$$

Observe that
$$
\sup_{u\in [0,A] } \frac1N \tilde Y_r(Nu) \to 0, \; a.s.
$$
for any $A\ge 0$. Since $\lambda_r$ are bounded, it follows that, for all $T>0$,
$$
\begin{array}{l}
\ds \sup_{t\in [0,T]} \left| \sum_{r\in\mR_{C}} \frac1N \gamma_{r}^C\tilde Y_{r}\left( N\int_{0}^t\tilde  \lambda_{r}(x^N_{C}(s))ds
\right) + \sum_{r\in S_{1}} \frac1N \gamma_{r}^C\tilde Y_{r}\left( N\int_{0}^t \tilde \lambda_{r}(x^N_{C}(s),X^N_{D}(s))ds\right)\right|\\
\\
 \to 0, \; a.s. \mbox{ when } N\to\infty.
\end{array}
$$
Clearly
$$
\sup_{t\in [0,T]} \left|\sum_{r\in \mR_{DC}\setminus S_{1}} \frac1N \gamma_{r}^CY_{r}\left( \int_{0}^t  \lambda_{r}(x^N_{C}(s),X^N_{D}(s))ds
\right) \right| \to 0, \; a.s.
$$
It follows that there exists a random constant $K_N$ going to zero such that, for $t, t_1,t_2\in [0,T]$, and  $\|F\|_\infty=\sup_{x\in \R^{M_C}\times \N^{M_D}}|F(x)|$.
$$
|x^N_{C}(t)|\le  |x^N_{C}(0)|+\|F\|_\infty t +K_N,\; a.s.$$
and
$$|x^N_{C}(t_1)-x^N_{C}(t_2)|\le \|F\|_\infty|t_1-t_2| +2K_N,\; a.s.$$

Tightness of $(x^N_C)_{N\in\N}$  in $C(\R^+; \R^{M_C})$ follows by classical criteria (see for instance \cite{jacod}, chapter 6, section 3b).

We conclude, from \cite{jacod} (chapter 6, section 3b), that $\{x^N\}_N=\{(x_C^N,X_D^N)\}_N$ is tight in $D(\R^+;E)$.

\par
\medskip
{\it{Step 2: Identification of limit points for bounded reaction rates.}}
\par

Let $x=(x_t	)_{t\geq 0}$ be the canonical process on $D(\R^+;E)$, and $P_{N}$ the law of $(x^N_t)_{t\geq 0}$ on this space, for each $N\in\N$.

We know that for each $N\in \N$ and $\varphi\in \mE.$

$$
\varphi (x_t) -\varphi(x_0)-\int_{0}^t \tilde \mA_{N } \varphi (x_s)ds
$$
is a $P_N$-martingale. Equivalently, for each $n\in \N$, $t_{1},\dots,t_{n}\in[0,r]$, $t\ge r \ge0$,
$\psi\in (C_{b} (E))^n$
and $\varphi\in \mE$
\begin{equation}
\label{e2}
\begin{array}{l}
\ds\E_{{P_N}}\left(\left( \varphi(x_{t})-\varphi(x_{0}) -\int_{0}^t \tilde \mA_{N } \varphi (x_{s})ds\right)
\psi(x_{t_{1}},\dots,x_{t_{n}})\right) \\
\ds=\E_{{P_N}}\left( \left(\varphi(x_{r})-\varphi(x_{0}) -\int_{0}^r \tilde \mA_{N } \varphi (x_{s})ds\right)
\psi(x_{t_{1}},\dots,x_{t_{n}})\right).
\end{array}
\end{equation}

Let $(P_{N_{k}})_{k}$ be a subsequence which converges weakly to  a measure $P$ on
$D(\R^+;E)$. We know that $x$ is $P$ almost surely continuous at every $t$ except
for a countable set $D_P$ and that for $t_{1}, \dots, t_{n}$ outside $D_P$, $P_{N_k}\pi^{-1}_{t_1,\ldots,t_n}$ converges weakly to $P\pi^{-1}_{t_1,\ldots,t_n}$ where $\pi_{t_1,\ldots,t_n}$ is the projection that carries the point $x\in D(\R^+;E)$ to the point $(x_{t_1},\ldots,x_{t_n})$ of $\R^n$.

Therefore, it is easy, using dominated convergence theorem and weak convergence properties, to let $k \to \infty$ in \eqref{e2} and obtain for  $t$, $t_{1},\dots, t_{n},r$ outside $D_P$:
 \begin{equation}
\label{e3}
\begin{array}{l}
\ds\E_{{P}}\left(\left( \varphi(x_{t})-\varphi(x_{0}) -\int_{0}^t  \mA_{\infty } \varphi (x_{s})ds\right)
\psi(x_{t_{1}},\dots,x_{t_{n}})\right) \\
\ds=\E_{{P}}\left( \left(\varphi(x_{r})-\varphi(x_{0}) -\int_{0}^r  \mA_{\infty } \varphi (x_{s})ds
\right)
\psi(x_{t_{1}},\dots,x_{t_{n}})\right).
\end{array}
\end{equation}

If $t\in D_P$, we choose a sequence $(t^k)$ outside $D_P$ such that $t^k\to t$ with $t^k>t$. Then $P\pi_{t^k}^{-1}$ converges weakly to $P\pi_{t}^{-1}$ since $x$ is $P$-a.s. right continuous in $t$ and $x_{t^k}$ converges almost surely to $x_t$. Then, we use \eqref{e3}
with $t^k$ instead of $t$, let $k\to  \infty$ and deduce that \eqref{e3} also holds for $t\in D_P$.
Similarly, we show that $t_{1},\dots, t_{n},r$ may be taken in $D_P$.

This shows that the measure $P$ is a solution of the martingale problem associated to the
generator $\mA_{\infty}$ on the domain $\mE$.

Hypothesis \ref{h2} enables us to apply Theorem \ref{t2.2}. The martingale problem has then a unique solution. It follows that the limit $P$ is equal to $P_{x_0}$, the law of the PDP, and that the whole sequence $(P_N)_N$ converges weakly to $P_{x_0}$.

\par
\medskip
{\it{Step 3: Conclusion}}
\par
Now, we prove Theorem \ref{t2.3.1} with a truncation argument.

Let $\theta \in C^\infty(\R^+)$ such that
$$
\left\{
\begin{array}{l}
\theta(x)=1,\; x\in [0,1],\\
\theta(x)=0,\; x\in [2,\infty),
\end{array}
\right.
$$
and, for $k\ge 1$ and $r\in\mR$, define
$$
\theta_k(x)=\theta\left(\frac{|x|^2}{k^2}\right),\; x\in E,
$$
and
$$
\lambda_r^k(x)=\theta_k(x)\lambda_r(x).
$$
Then, the problem with $\lambda_r^k$ instead of $\lambda_r$ fulfills Hypothesis \ref{h2}. We define $x^N_k=(x^N_{C,k},X^N_{D,k})$ the jump Markov process associated to the jump intensities $\lambda_r^k$, starting at $x^N(0)$. By the preceding result, we know that, for all $k\in\N$, $(x^N_{k})_{N\in\N}$ converges weakly to the PDP $x_{k}$ in $D(\R^+;E)$, whose characteristics are the jump intensities $\lambda_r^k$, with corresponding transition measure, and vector fields (obvious definitions as in Theorem \ref{t2.3.1}).

Then, $((x^N_{k})_{k\in\N})_{N\in\N}$ converges weakly to $(x_{k})_{k\in\N}$ in $D(\R^+;E)^\N$.

By Skorohod representation Theorem (see \cite{bill2}, Theorem 3.3), up to a change of  probability space, we may assume that that for all $k\in\N$
$(x^N_{k})_N$ converges a.s. to
$x_{k}$ in $D(\R^+;E)$.

Let $T>0$ and the stopping times
$$
\tau^k=\inf\{ t\in[ 0,T],\; |x_k(t)|\ge k\},
$$
with $\tau^k=T$ if $\{ t\in[ 0,T],\; |x_k(t)|\ge k\}=\emptyset$.

Then, for $k,l\in\N$
$$
x_{k}(t)=x_{l}(t), \; t\in [0,\tau^l\wedge \tau^k], \; a.s.
$$
so that $\tau^k$ is a.s. non-decreasing.

Moreover, if $x$ (resp. $x^N$) are the PDP associated to $\mathcal A_\infty$ (resp. the Markov jump process associated to $\tA_N$), then
$$
x_{k}(t)=x_t,\; t\in [0,\tau^k), \; a.s.
$$
and if
$$
\tau^k_N=\inf\{ t\in[ 0,T],\; |x_k^N(t)|\ge k\},$$ with $\tau^k_N=T$ if $\{ t\in[ 0,T],\; |x_k^N(t)|\ge k\}=\emptyset$, then
$$
x_k^N(t)=x^N(t),\; t\in [0,\tau^k_N), \; a.s.
$$

Let $\delta>0$. Observing that if $\tau^{k-1}>T-\delta$ and $d_{T-\delta}(x^N_k,x_k)<\epsilon$, where $d_{T-\delta}$ is the distance on $D([0,T-\delta];E)$, then, for enough small $\epsilon$:
$$
\sup_{t\in [0,T-\delta]}  |x_{k}^N(t)| \le k, \; a.s.
$$
then a.s., $\tau^k_N\geq T-\delta$ and $x^N_{k}=x^N$ in $[0,T-\delta]$. Since $\tau^k\ge \tau^{k-1}>T-\delta$; we have also, a.s.,
$x_{k}=x$ in $[0,T-\delta]$ and
$$
d_{T-\delta}(x^N,x)<\epsilon.
$$
We deduce that
$$
\forall \delta>0, \quad \P\left( d_{T-\delta}(x^N,x)\ge \epsilon \right)\le \P\left( \tau^{k-1}\le T-\delta \right)
+\P\left( d_{T-\delta}(x^N_k,x_k)\ge \epsilon \right).
$$
Finally, we have
$$
\P\left( \tau^{k-1}\le T-\delta \right) = \P(\sup_{t\in [0,\tau^k]}|x_k(t)|\ge k-1)\le
\P(\sup_{t\in [0,T-\delta]}|x(t)|\ge k-1).
$$
By Hypothesis \ref{h1}, the PDP $x$  can not explode in finite time. Thus for   $k$ large, this term
is small. Then for large $N$, the second is small. We deduce that $x^N$ converges in probability to $x$ in the new probability space. Returning in the original probability space, we obtain that $x^N$ converges in distribution to $x$.

\section{Piecewise deterministic process with jumps}

In this section, we assume that some of the reactions in a subset $S_{2}$ of $\mR_{DC}$ are such that
$\gamma_{r}^C$ is large and scales with $N$. We set $\tilde \gamma_{r}^C= \frac1N\gamma_{r}^C$ for
$r\in S_{2}$. We define $S=S_1\cup S_2$.

The scaled generator has now the form
\begin{equation}
\label{e4.1}
\begin{array}{ll}
\tilde \mA_{N} f(x_{C},X_{D}) &=\ds  \sum_{r\in \mR_{C}} \left[ f(x_{C}+\frac1N \gamma_r^C,X_{D})-f(x_{C},X_{D})
\right]N\tilde\lambda_r(x_{C})\\
\\
&\ds +\sum_{r\in S_{2}} \left[ f(x_{C}+\tilde \gamma_r^C, X_{D}+\gamma_{r}^D)-f(x_{C},X_{D})\right] \lambda_r(x_{C},X_{D})\\
\\
&\ds +\sum_{r\in S_{1}} \left[ f(x_{C}+\frac1N \gamma_r^C, X_{D})-f(x_{C},X_{D})\right] N\tilde\lambda_r(x_{C},X_{D})\\
\\
&\ds +\sum_{r\in \mR_{DC}\setminus S} \left[ f(x_{C}+\frac1N\gamma_r^C,X_{D}+\gamma_{r}^D)-f(x_{C},X_{D})\right]\lambda_r(x_{C},X_{D})\\
\\
&\ds +\sum_{r\in \mR_{D}} \left[ f(x_{C}, X_{D}+\gamma_r^D)-f(x_{C},X_{D})\right]\lambda_r(X_{D}).
\end{array}
\end{equation}

For $f\in C^1_{b}(E)$, we may let $N\to \infty$ and obtain the limit
generator
$$
\begin{array}{ll}
\ds \tA_{\infty}f(x_{c},X_{D}) &=\ds  \left(\sum_{r\in \mR_{C}} \tilde \lambda_{r}(x_{C})  \gamma_{r}^C
+\sum_{r\in S_1} \tilde \lambda_{r}(x_{C},X_D)  \gamma_{r}^C\right)
\cdot \nabla_{x_{C}}f(x_{C},X_{D}) \\
\\
&\ds + \sum_{r\in S_{2}}  \left[ f(x_{C}+\tilde \gamma_r^C, X_{D}+\gamma_{r}^D)-f(x_{C},X_{D})\right] \lambda_r(x_{C},X_{D})\\
\\
&\ds  +\sum_{r\in \mR_{DC}\setminus S} \left[ f(x_{C},X_{D}+\gamma_{r}^D)-f(x_{C},X_{D})\right]\lambda_r(x_{C},X_{D})\\
\\
&\ds +\sum_{r\in \mR_{D}} \left[ f(x_{C}, X_{D}+\gamma_r^D)-f(x_{C},X_{D})\right]\lambda_r(X_{D}).
\end{array}
$$
This formal argument indicates that if $N\to\infty$, the process is a  piecewise
deterministic process with jumps both in $x_C$ and $X_D$. In fact, the proof of this can be done easily thanks to the result of section
3.

Indeed, let us introduce the following auxiliary system of reactions involving the variable $(x^{1,N}_{C},x^{2,N}_{C},X^N_D)\in\R^{M_C}\times \R^{M_C}\times \N^{M_D}$:
\begin{itemize}
\item If $r\in \mR_{C}$, $\ds (x^{1,N}_{C},x^{2,N}_{C},X_D^N)\to (x^{1,N}_{C}+\frac1N \gamma^{r}_{C},x^{2,N}_{C},X_D^N)$\\
\item If $r\in S_{2}$, $(x^{1,N}_{C},x^{2,N}_{C},X_D^N)\to (x^{1,N}_{C},x^{2,N}_{C}+\tilde \gamma^{r}_{C},X_D^N+\gamma^r_D)$\\
\item If $r\in \mR_{DC}\setminus S_2$, $\ds (x^{1,N}_{C},x^{2,N}_{C},X_D^N)\to (x^{1,N}_{C}+\frac1N \gamma^{r}_{C},x^{2,N}_{C},X_D^N+\gamma^r_{D})$\\
\item If $r\in \mR_{D}$, $\ds (x^{1,N}_{C},x^{2,N}_{C},X_D^N)\to (x^{1,N}_{C},x^{2,N}_{C},X_D^N+\gamma^r_{D})$
\end{itemize}
\par
\medskip
The rates of these reactions are $\lambda_r(x^{1,N}_{C},x^{2,N}_{C},X_D^N)=\lambda_r(x^{1,N}_{C}+x^{2,N}_{C},X_D^N)$.

If we choose the initial data $(x_{C}^N(0),0,X^N_{D}(0))$ then this system and the original
system are equivalent. Indeed, we recover the value of the original system by the
relation:
$$
(x^{N}_{C},X_D^N)=(x^{1,N}_{C}+x^{2,N}_{C},X^N_D).
$$
Conversely, notice that $x^{N}_{C}$ is a sum of a pure jump part and of a continuous one so that   given $(x^{N}_{C},X_D^N)$ , $(x^{1,N}_{C},x^{2,N}_{C},X_D^N)$ can be recovered
by isolating from $x^N_C$ the small jumps to obtain $x^{1,N}_{C}$ and the jump of order one
(from $S_2$) to obtain $x^{2,N}_{C}$.

The auxiliary system corresponds to the following generator:
\begin{equation}
\label{e1bis}
\begin{array}{l}
\mA_{aux}f(x^{1}_{C},x^{2}_{C},X^D) \\
\\
=\ds  \sum_{r\in \mR_{C}} \left[ f(x^1_{C}+\frac1N \gamma_r^C,x^2_C,X_{D})-f(x^1_{C},x^2_C,X_{D})
\right]N\tilde\lambda_r(x^1_{C}+x^2_C)\\
\\
\ds +\sum_{r\in S_{2}} \left[ f(x^1_{C},x^2_C+ \tilde \gamma_r^C, X_{D}+\gamma_{r}^D)-f(x^1_{C},x^2_C,X_{D})\right] \lambda_r(x^1_{C}+x^2_C,X_{D})\\
\\
\ds +\sum_{r\in S_{1}} \left[ f(x^1_{C}+\frac1N \gamma_r^C,x^2_C, X_{D})-f(x^1_{C},x^2_C,X_{D})\right] N\tilde \lambda_r(x^1_{C}+x^2_C,X_{D})\\
\\
\ds +\sum_{r\in \mR_{DC}\setminus S} \left[ f(x^1_{C}+\frac1N\gamma_r^C,x^2_C,X_{D}+\gamma_{r}^D)-f(x^1_{C},x^2_C,X_{D})\right]\lambda_r(x^1_{C}+x^2_C,X_{D})\\
\\
\ds +\sum_{r\in \mR_{D}} \left[ f(x^1_{C},x^2_C, X_{D}+\gamma_r^D)-f(x^1_{C},x^2_C,X_{D})\right]\lambda_r(X_{D}).
\end{array}
\end{equation}
This generator is of the form as in section 3 but the discrete variable is now $(x^2_C,X_D)$.
Using the results proved in this section we thus obtain that this auxiliary system of reaction
converges in distribution to the piecewise continuous deterministic process given by the generator
$$
\begin{array}{l}
\ds \mA_{aux,\infty}f(x_{C}^{1},x_{C}^{2},X_D) \\
=\ds  \left(\sum_{r\in \mR_{C}} \tilde \lambda_{r}(x^1_{C}+x^2_C,X_D)  \gamma_{r}^C
+\sum_{r\in S_1} \tilde \lambda_{r}(x^1_{C}+x^2_C,X_D)  \gamma_{r}^C\right)
\cdot \nabla_{x^1_{C}}f(x^1_{C},x^2_C,X_{D}) \\
\\
\ds + \sum_{r\in S_{2}}  \left[ f(x^1_{C},x^2_C+\tilde \gamma_r^C, X_{D}+\gamma_{r}^D)-f(x^1_{C},x^2_C,X_{D})\right] \lambda_r(x^1_{C}+x^2_C,X_{D})\\
\\
\ds  +\sum_{r\in \mR_{DC}\setminus S} \left[ f(x^1_{C},x^2_C,X_{D}+\gamma_{r}^D)
-f(x^1_{C},x^2_C,X_{D})\right]\lambda_r(x^1_{C}+x^2_C,X_{D})\\
\\
\ds +\sum_{r\in \mR_{D}} \left[ f(x^1_{C},x^2_C, X_{D}+\gamma_r^D)-f(x^1_{C},x^2_C,X_{D})\right]\lambda_r(X_{D}).
\end{array}
$$
Going back to the original variables, {i.e.} setting $x_C=x^1_{C}+x^2_C$, we deduce that the original system converges in distributions to the piecewise deterministic process with generator
$\tilde\mA_\infty$.

To obtain this convergence, we suppose that the jump intensities are  $C^1$-differentiable with respect to the variable $x_C^1$ and that the limit PDP satisfies Hypothesis
\ref{h1}.

% We suppose also that the jump intensities and their derivatives with respect to $x_C^1$ are bounded on $\R^{M_C}\times\R^{M_C}\times\N^{M_D}$.

\section{Averaging}
In this section, we examine the case when  (\ref{e2.3.1}) is not satisfied by all the discrete variables. In this case, the previous results are not valid. We introduce the decomposition $X_D=(X_D^1,X_D^2)\in \mathbb N^{M_{D,1}}\times\mathbb N^{M_{D,2}}$ and $\gamma_r^D=(\gamma_r^{D,1},\gamma_r^{D,2})$ such that
$$\gamma_r^{D,1}=0,\qquad r\in S_1.$$
This replaces (\ref{e2.3.1}). The continuous variable $x_C$ follows the same characteristics as in section \ref{CPDP}.

In the set of the reactions $S_1$, the discrete variable $X_D^2$ has fast motion, its jumps rates are of order $N$ and its jumps are of order 1.

More precisely, for $r\in S_1$, the state changes according to $$(x_C,X_D^1,X_D^2)\to (x_C+\frac1N\gamma_r^C,X_D^1,X_D^2+\gamma_r^{D,2})$$ with the rate $N\widetilde\lambda_r(x_C,X_D^1,X_D^2)$.

We now have the following  generator for the process :

\begin{equation}
\label{e2.4.1}
\begin{array}{l}
\tA_{N} f(x_{C},X_{D}^1,X_D^2) \\
\\
=\ds  \sum_{r\in \mR_{C}} \left[ f(x_{C}+\frac1N \gamma_r^C,X_{D}^1,X_D^2)-f(x_{C},X_{D}^1,X_D^2)
\right]N\tilde\lambda_r(x_{C})\\
\\
\ds +\sum_{r\in S_{1}} \left[ f(x_{C}+\frac1N\gamma_r^C,X_{D}^1,X_D^2+\gamma_r^{D,2})-
f(x_{C},X_{D}^1,X_D^2)\right]N\tilde\lambda_r(x_{C},X_{D}^1,X_D^2)\\
\\
\ds +\sum_{r\in \mR_{DC}\setminus S_{1}}
\left[ f(x_{C}+\frac1N\gamma_r^C,X_{D}^1+\gamma_r^{D,1},X_D^2+\gamma_r^{D,2})
-f(x_{C},X_{D}^1,X_D^2)\right]\lambda_r(x_{C},X_{D}^1,X_D^2)\\
\\
\ds +\sum_{r\in \mR_{D}} \left[ f(x_{C}, X_{D}^1+\gamma_r^{D,1},X_D^2+\gamma_r^{D,2})
-f(x_{C},X_{D}^1,X_D^2)\right]\lambda_r(X_{D}^1,X_D^2).
\end{array}
\end{equation}

This new model combines slow and fast motions. This leads to double time scale evolution which can be further simplified.
%is extremely difficult to analyse directly.
Contrary to the previous sections, we can not take the formal limit of this generator.
The idea is to average in fast discrete variables $X_D^2$, and focus on the slow variables ($x_C,X_D^1)$, in order to obtain a much simpler averaged generator. We assume that $X_D^2$ takes only a finite number of values: there exists a finite set $K$ in $\N^{M_{D,2}}$ such that
$X_D^2\in K$.

To this aim, we introduce the following generator depending on $(x_C,X_D^1)$:
$$
\mA_{x_C,X_D^1}h(X_D^2)= \sum_{r\in S_1} \left[ h(X_{D}^2+\gamma_r^{D,2})
-h(X_{D}^2)\right]\tilde \lambda_r(x_C,X_{D}^1,X_D^2).
$$

We assume that for all $(x_C,X_D^1)$, the process associated to the generator $\mA_{x_C,X_D^1}$ is uniquely ergodic, and we denote by $\nu_{x_C,X_D^1}$ its unique invariant measure.

Then we define the averaged jump rates:
\begin{equation}
\label{e2.4.5}
\begin{array}{l}
\ds\bar\lambda_r(x_{C},X_{D}^1) =\int_{\N^{M_{D,2}}} \tilde\lambda_r(x_{C},X_{D}^1,X_D^2)
\nu_{x_C,X_D^1}(dX_D^2),
\; r \in  S_1,\\
\ds\bar\lambda_r(x_{C},X_{D}^1) =\int_{\N^{M_{D,2}}} \lambda_r(x_{C},X_{D}^1,X_D^2)
\nu_{x_C,X_D^1}(dX_D^2),
\; r \in \mR_{DC}\setminus S_1,\\
\ds\bar\lambda_r(x_C,X_{D}^1) =\int_{\N^{M_{D,2}}} \lambda_r(X_{D}^1,X_D^2)
\nu_{x_C,X_D^1}(dX_D^2),
\; r \in \mR_D.
\end{array}
\end{equation}

We present now the main result of this section.
\begin{Theoreme}
\label{t2.4.1}
Let $x^N=(x_C^N,X^N_{D,1},X^N_{D,2})$ be a jump Markov process with generator $\tA_N$, starting at $x^N(0)=(x_C^N(0),X_{D,1}^{N}(0),X_{D,2}^{N}(0))$. Assume that the set of values of $X_{D,2}$
is finite and that for all $(x_C,X_D^1)$, the process associated to the generator $\mA_{x_C,X_D^1}$ is
uniquely ergodic. Assume also that the jump rate $\tilde \lambda_r$, $r\in \mR_C\cap S_1$,
and $\lambda_r$, $r\in \mR_{DC}\setminus S_1$ and $\bar\lambda_r$, $r\in \mR_{DC}\cup \mR_D$, are $C^1$ with respect to $x_C$.
 Define $P_{x_0}$ the law of the PDP starting at $x_0=(x_{C,0},X_{D,0}^1)$ whose jump intensities are:
$$\lambda(x_C,X_D^1)=\sum_{r\in (\mR_{DC}\setminus S_1)\cup \mR_D}\bar{\lambda}_r(x_C,X_D^1),$$
with transition measure:
$$\int_{E}f(z)Q(dz;x_C,X_D^1)=\frac{1}{\lambda(x_C,X_D^1)}\sum_{r\in (\mR_{DC}\setminus S_1)\cup\mR_D }f(x_C,X_D^1+\gamma_r^{D,1})\bar{\lambda}_r(x_C,X_{D}^1),$$
and vectors fields:
$$F_{X_{D}^1}(x_C)=\sum_{r\in\mR_C}\gamma_r^C\tilde\lambda_r(x_C)+\sum_{r\in S_1}\gamma_r^C\bar\lambda_r(x_C,X_D^1).$$

Assume  that
\begin{enumerate}
\item $\lambda, Q$ and $F$ satisfy Hypothesis \ref{h1}.%, and the differential system:
% $$
% \left\{
%\begin{array}{l}
%\ds\frac{dx_C}{dt}= F_{X_D^1}(x_C),\\
%x_C(0)=x_{C,0},
%\end{array}
%\right.
%$$
%has a unique global solution, for every $X^1_D$.
\item $(x_C^N(0),X_D^{1,N}(0))$ converges in distribution to $x_0$ in $\R^{M_C}\times\N^{M_{D,1}}$.
\item there exists $K_1>0$ and $K_2>0$ such that, for every $s\in\R^+,(x_C,X_D^1,X_D^2)\in E$, and for every bounded function $g$ satisfying the centering condition \begin{equation}\label{CenteringCond}
\int_{\N^{M_{D,2}}} g(X_D^2)\nu_{x_C,X_D^1}(dX_D^2)=0\end{equation} we have $$P_s^{x_C,X_D^1}g(X_D^2)<K_1e^{-K_2s}\| g\|_\infty$$ where $(P_s^{x_C,X_D^1})_s$ is the semigroup associated to $\mA_{x_C,X_D^1}$.
\item if $g$ is a bounded function satisfying the centering condition (\ref{CenteringCond}), the Poisson equation
$$
\sum_{r\in S_{1}} \left[ h(x_{C},X_{D}^1,X_D^2+\gamma_r^{D,2})-
h(x_{C},X_{D}^1,X_D^2)\right]\tilde\lambda_r(x_{C},X_{D}^1,X_D^2)=g(X_D^2)
$$
has a solution given by
$$h(x_C,X_D^1,X_D^2)=-\int_0^\infty P_s^{x_C,X_D^1}g(X_D^2)ds$$ and this solution is Lipschitz with respect to $x_C$, uniformly in $X_D^1$.

\end{enumerate}

Then $(x^N_C,X_{D}^{1,N})$ converges in distribution in $D(\mathbb R^+;\mathbb R^{M_C}\times\mathbb N^{M_{D,1}})$ to the PDP whose law is $P_{x_0}$.
\end{Theoreme}
\begin{Remark}
Assumptions 1. and 2. are similar to those made previously. Assumptions 3. and 4. guaranty that it is
possible to average with respect to the variable $X^2_D$. Assumption 3. states that the dynamic associated to the generator $\mA_{x_C,X_D^1}$ is exponentially mixing, uniformly with respect to
$(x_C,X^1_D)$, assumption 4. then follows if $\mA_{x_C,X_D^1}$ depends smoothly on $x_C$.
\end{Remark}
{\bf{Proof}}

The steps of the proof are similar to the proof of the Theorem \ref{t2.3.1}. We begin by supposing that all the jump rates are bounded as well as their derivatives. This Hypothesis will be removed at the end of the proof. We take
$E=\R^{M_C}\times \N^{M_{D,1}}$.

{\it{Step 1: Tightness}}.
\par

The proof of the tightness of the laws of $(x_C^N,X_{D,1}^{N})$ is similar as in section \ref{CPDP}, and is left to the reader.

{\it{Step 2: Identification of limit points}}.
\par
Since we are interested by the evolution of ($x_C,X_D^{1})$, it is natural to consider, for the generator $\tA_N$, test functions $f(x_C,X^1_D)$ depending only on those variables. Unfortunately, it remains some terms of order 1 depending on the variable $X_D^2$. To overcome this problem, we introduce 'perturbed test functions' (see \cite{papanicolaou-strook-varadhan}, \cite{kushner}, \cite{fouque}) :
$$
f_N(x_C,X_D^1,X_D^2)=f(x_C,X_D^1)+\frac1N f^1(x_C,X_D^1,X_D^2),
$$
with $f\in\mathcal E$ and $f$ independent of $X_D^2$. We temporary need
that $f$ is also $C^2$ with bounded second derivatives with respect to $x_C$. The function $f^1$ is chosen such that the formal limit of $\tA_N f_N(x_C,X_D^1,X_D^2)$ does not depend on $X_D^2$. 

This can be down by taking $f^1$ such that
\begin{equation}
\label{e2.4.4}
\sum_{r\in S_{1}} \left[ f^1(x_{C},X_{D}^1,X_D^2+\gamma_r^{D,2})-
f^1(x_{C},X_{D}^1,X_D^2)\right]\tilde\lambda_r(x_{C},X_{D}^1,X_D^2)=g_{x_{C},X_{D}^1}(X_D^2)
\end{equation}

with \begin{equation}
\label{e2.4.4bis}
\begin{array}{l}
g_{x_C,X_D^1}(X_D^2) \\
\\
=\ds
\left(\sum_{r\in S_{1}} \gamma_r^C \bar \lambda_r(x_{C},X_{D}^1)- \gamma_r^C\tilde\lambda_r(x_{C},X_{D}^1,X_D^2)
\right)  \cdot\nabla_{x_C} f(x_C,X_D^1)\\
\\
\ds +\sum_{r\in \mR_{DC}\setminus S_{1}}
\left[ f(x_{C},X_{D}^1+\gamma_r^{D,1})
-f(x_{C},X_{D}^1)\right]\left( \bar \lambda_r(x_{C},X_{D}^1)-\lambda_r(x_{C},X_{D}^1,X_D^2)\right)\\
\\
\ds +\sum_{r\in \mR_{D}} \left[ f(x_{C}, X_{D}^1+\gamma_r^{D,1})
-f(x_{C},X_{D}^1)\right]\left( \bar \lambda_r(x_C,X_{D}^1)-\lambda_r(X_{D}^1,X_D^2)\right)
\end{array}
\end{equation}
Since the jump rates are assumed to be bounded, $g_{x_C,X_D^1}$ is clearly bounded. Also,
the function $g_{x_C,X_D^1}$ is zero mean with respect to the invariant distribution, i.e. $g_{x_C,X_D^1}$ satisfy (\ref{CenteringCond}).

%For all $(x_C,X_D^1)$, the infinitesimal generator $\mA_{x_C,X_D^1}$ is not invertible, since it has the non trivial function $h(X_D^2)=1$ in its null space, that is $\mA_{x_C,X_D^1} 1=0$. Moreover, the null space of $\mA_{x_C,X_D^1}$ is reduced to the constant functions, by ergodicity of the Markov Process associated to $\mA_{x_C,X_D^1}$.

Then, the function $f^1$ satisfying \eqref{e2.4.4} exists under our assumptions. Moreover, it is not difficult to check that $f^1$ is Lipschitz with respect to $x_C$, uniformly with respect to $X_D^1$.

\medskip
{\it{Limiting infinitesimal generator}.}
\noindent

\par
\medskip

We now return to the analysis of the limiting problem.
 We denote by $P_{x,N}$ the law of $(x_C^N,X_{D,1}^{N}, X_{D,2}^{N})$, then
$$
\begin{array}{c}
\ds f_N(x_C(t),X_D^1(t),X_D^2(t))-f_N(x_C(0),X_D^1(0),X_D^2(0)) \\
\\
\ds -\int_0^t \tA_Nf_N(x_C(s),X_D^1(s),X_D^2(s))ds
\end{array}
$$
is a $P_{x,N}$ martingale.
Equivalently, for each $n\in \N$, $t_{1},\dots,t_{n}\in[0,r]$, $t\ge r \ge0$,
$\psi\in (C_{b} (E))^n$

\begin{equation}
\label{e2bis}
\begin{array}{l}
\ds\E_{{P_{x,N}}}\left(\left( f_N(x_{t})-f_N(x_{0}) -\int_{0}^t \tilde \mA_{N } f_N (x_{s})ds\right)
\psi(x_{t_{1}},\dots,x_{t_{n}})\right) \\
\ds=\E_{{P_{x,N}}}\left(\left( f_N(x_{r})-f_N(x_{0}) -\int_{0}^r \tilde \mA_{N } f_N (x_{s})ds\right)
\psi(x_{t_{1}},\dots,x_{t_{n}})\right)
\end{array}
\end{equation}

Let us now consider  a subsequence $(\tilde P_{x,N_k})_k$ of $(\tilde P_{x,N})_N$, the laws of ($x_C^N,X_{D,1}^{N})$ which weakly converges to a measure $\tilde P_x$ on $D(\mathbb R^+;\mathbb R^{M_C}\times \mathbb N^{M_{D,1}})$.

%Moreover, it can be deduced from Lemma \ref{lem} that $f^1$ is $C^1$ with respect to $x_C$ and its gradient with respect to $x_C$ is uniformly bounded on $E$. Also, $f\in \mathcal E$, $f$ is $C^2$ with respect to $x_C$, and its second derivative with respect to $x_C$ is uniformly bouned.

%
%
%
%\begin{Lemma}
%\label{lem}
%The function $(x_C,X_D^1,X_D^2)\mapsto \nu_{x_C,X_D^1}(X_D^2)$ is bounded, $C^1$ differentiable with respect to the continuous variable $x_C$ and its derivatives are uniformly bounded in $E$.
%\end{Lemma}
%
%\textbf{Proof of the Lemma:}
%\noindent
%The boundness is trivial as the functions $\lambda_r$ are supposed bounded. By ergodicity, $\nu_{x_C,X_D^1}$ is a simple eigenvector of  $\mA_{x_C,X_D^1} $ for the eigenvalue 0. The implicit function theorem ensures the existence of an eigenvalue function $(x,y)\mapsto \lambda_1(x,y)$ defined on a neighborhood of $x_C,X_D^1$ and $C^1$-differentiable with respect to $x$ and $y$. Moreover $\lambda_1(x,y)$ is a simple eigenvalue of $\mA_{x,y}$. The corresponding eigenvector $v_1(\mA_{x,y})$ can be calculated by fixing one of its component to 1, and resolving a Cramer system. The components of $v_1(A)$ are then products of the principal determiner of $\mA_{x,y}$ and of components of $\mA_{x,y}$. We can deduce that $v_1$ is $C^1$ with respect to the components of $x,y$.
%
%$\square$

It is clear that $f_N$ converges to $f$, uniformly on $E$. Moreover, the continuity and boundedness properties of $f^1$ ensure that $\tA_Nf_N$ converges on $E$ towards $\tA_\infty f$, where

\begin{equation}
\label{e2.4.6}
\begin{array}{l}
\tA_{\infty} f(x_{C},X_{D}^1) \\
\\
=\ds  \left( \sum_{r\in \mR_{C}}\gamma_r^C \tilde\lambda_r(x_{C}) +
\sum_{r\in S_{1}} \gamma_r^C\bar\lambda_r(x_{C},X_{D}^1)\right) \cdot\nabla_{x_C} f(x_C,X_D^1)\\
\\
\ds +\sum_{r\in (\mR_{DC}\setminus S_1)\cup \mR_{D}} \left[ f(x_{C}, X_{D}^1+\gamma_r^{D,1})
-f(x_{C},X_{D}^1)\right]\bar\lambda_r(x_C,X_{D}^1).
\end{array}
\end{equation}
It is the generator of a PDP whose characteristics have been averaged in fast variable $X_D^2$.

Equation (\ref{e2bis}) gives, as $N_k\to\infty$:
\begin{equation}
\label{17}
\begin{array}{l}
\ds\E_{{\tilde{P}_x}}\left(\left( f(x_{t})-f(x_{0}) -\int_0^t \tilde \mA_\infty f (x_{s})ds\right)
\psi(x_{t_{1}},\dots,x_{t_{n}})\right) \\
\ds=\E_{{\tilde P_x}}\left(\left( f(x_{r})-f(x_{0}) -\int_{0}^r \tilde \mA_{\infty } f (x_{s})ds
\right)
\psi(x_{t_{1}},\dots,x_{t_{n}})\right).
\end{array}
\end{equation}
Recall that we have assumed that $f$ is $C^2$ with respect to $x_C$. It is easy to prove 
that \eqref{17}
holds for any $f\in \mE$ by choosing a sequence $(f_n)$in $\mE$ such that, for all $n$, $f_n$
is $C^2$ with bounded second derivatives and $f_n$ and $\nabla f_n$ are uniformly bounded
and
converge pointwise to $f$ and $\nabla f$  . 

The measure $\tilde P_x$ is then a solution of the martingale problem associated to the generator $\tilde\mA_{\infty}$ on the domain of the functions in $\mE$ independent of $X_D^2$.

By the Theorem \ref{t2.2}, the martingale problem has then a unique solution. We can then deduce that $(\tilde P_{x,N})_N$ converges weakly to $P_{x_0}$, the law of the PDP defined in the Theorem \ref{t2.4.1}.

{\it{Step 3: Conclusion}}.
\par
Since $X_D^2$ is assumed to live in a finite set, we can remove the boundedness assumption
on the jumps rates as in section 3. We then conclude as in section \ref{CPDP}. Theorem \ref{t2.4.1} is proved.

\section{Discontinuous PDP and singular switching}

In this section, our system has two time scales. The switching between fast and slow dynamics
is governed by the state of a discrete variable. For simplicity, we suppose now that there is no
other discrete variable and denote by $\theta$, taking values in $\{0,1\}$, this unique discrete
variable which governs the time scale of the continuous variables.

The number of molecules belonging to species $C$ are again
supposed to be of order $N$, and we continue to write $x_C=\frac{X_C}{N}$.

When $\theta=0$, the rates of all reactions from $\mathcal R_{C} \cup \mathcal R_{DC}$ is of order $N$. When $\theta=1$, some reactions (inactive for $\theta=0$)
become active with much faster rates of order $\frac{N}{\epsilon}$, where $\epsilon$ is a new parameter, supposed to be small.
For simplicity we consider that there are no $\mathcal R_{C}$ type reactions. This would
cause only slightly more complicated notations.  Then for $r\in\mathcal R_{DC}$, we introduce
the rescaled reaction rates
$$
\begin{array}{ll}
\ds \tilde \lambda_r(x_C,0)=&\ds \frac1N\lambda_r(x_C,0)\\
\\
\ds \tilde \lambda_r(x_C,1)=&\ds \frac{\epsilon}{N}\lambda_r(x_C,1).
\end{array}$$

%\corrO{
%These rates satisfy, for $r \in S_3$
%$$
%\begin{array}{ll}
%\tilde \lambda_r(x_C,0) = 0, &\tilde \lambda_r(x_C,1) > 0.
%\tilde \lambda_r(x_C,0) > 0, &\tilde \lambda_r(x_C,1) = 0, &\text{if}\quad r \in \mathcal R_{DC \& C} \setminus S_3.
%\end{array}
%$$}

This system will not have a limit if $\theta$ stays for too long time in the state 1. The rate of the reaction changing $\theta$ from state 1 to state 0 is supposed to be of order $\frac1\epsilon$ and
we also set:
$$ \tilde\lambda_\theta(x_C,1)=\epsilon\lambda_{\theta}(x_C,1).$$
The opposite reaction changing $\theta$ from state 0 to state 1 is written $\lambda_\theta(x_C,0)$.

In these conditions, the generator of the process has the following form:
$$
\begin{array}{l}
\tilde{\mathcal A}_{N,\epsilon}  f(x_{C},0) \\
\ds =  \sum_{r\in {R_{DC}}} \left[ f(x_{C}+\frac1N \gamma_r^C,0)-f(x_{C},0)
\right]N\tilde\lambda_r(x_{C},0)
+  \left[ f(x_{C},1)-f(x_{C},0)
\right]\lambda_\theta(x_{C},0),
\end{array}
$$
and
$$
\begin{array}{l}
\tilde {\mathcal A}_{N,\epsilon}  f(x_{C},1)\\
\ds = \frac1\epsilon \sum_{r\in {\mathcal R_{DC}}} \left[ f(x_{C}+\frac1N \gamma_r^C,1)-f(x_{C},1)
\right]N\tilde\lambda_r(x_{C},1)+   \frac1\epsilon \left[ f(x_{C},0)-f(x_{C},1)
\right]\tilde\lambda_\theta(x_{C},1).
\end{array}
$$

At the limit (for high $N$ and small $\epsilon$), we will show that the discrete process $\theta$ inducing kicks in the continuous variable $x_C$ is almost surely equal to 0.
We are interested in the limit distribution of the process $x_C$.

\par
\medskip

We introduce the flow $\phi_1(t,x_C)$, associated to the vector fields $$F_1(x_C)=\sum_{r\in R_{DC } }\gamma_r^C\tilde{\lambda}_r(x_C,1).$$

We now state the main result of this section.
\begin{Theoreme}
Let $x^{N,\epsilon}=(x_C^{N,\epsilon},\theta^{N,\epsilon})$ be a jump Markov process with values in $\R^{M_C}\times\{0,1\}$, with generator $\tA_{N,\epsilon}$, starting at $x^{N,\epsilon}(0)=
\left(x_C^{N,\epsilon}(0),\theta^{N,\epsilon}(0)\right)$.
Assume that the following assumptions hold:
\begin{itemize}
\item $x_C^{N,\epsilon}(0)$ converges in distribution to $x_C(0)$ in $\R^{M_C}$ as $N\to\infty$
and $\epsilon\to0$.
\item $\theta^{N,\epsilon}(0)$ converges in distribution to $0$ as $N\to\infty$ and $\epsilon\to0$.
\item The jump rates $\tilde \lambda_r, $ for $r\in \mathcal R_{DC }$, $\tilde \lambda_\theta(\cdot,1)$ and $\lambda_\theta(\cdot,0)$ are uniformly bounded, $C^1$ with respect to $x_C$, and their derivatives are uniformly bounded. Moreover, there exists $\alpha>0$, such that
$$\tilde\lambda_\theta(x_C,1)\geq \alpha,\qquad \forall x_C\in\R^{M_C}.$$
\end{itemize}

Then $(x^{N,\epsilon}_C)$ converges in distribution in $L^p([0,T];\mathbb R^{M_C})$,
for any $\infty>p\ge 1$, to the PDP whose generator is given by
$$
\begin{array}{l}
\tA_{\infty}  \varphi(x_{C})
\ds = \left( \sum_{r\in \mathcal R_{DC } }\gamma_r^C\tilde\lambda_r(x_{C},0) \right)\cdot\nabla_{x_C}\varphi(x_{C})\\
\\
\ds +\lambda_\theta(x_{C},0) \int_0^\infty \left(\varphi(\phi_1(t,x_C))-\varphi(x_{C})\right)\tilde\lambda_\theta(\phi_1(t,x_C),1)e^{-\int_0^t
\tilde\lambda_\theta(\phi_1(s,x_C),1)ds} dt\\
\end{array}
$$ for every $\varphi\in C_b^1(\R^{M_C})$.
\end{Theoreme}

\begin{Remark}

%$\tA_{\infty}$ is the generator of a PDP  whose characteristics are
%$\sum_{r\in \mathcal R_{DC } }\gamma_r^C\tilde\lambda_r(x_{C},0)$, $\lambda_\theta(\cdot,0)$, and $Q$ defined by $Q(dz; x)=\tilde\lambda_\theta\left(\phi_1(t,x),1\right)e^{-\int_0^t\tilde\lambda_\theta\left(\phi_1(s,x),1\right)ds}dt$ with $z=\phi_1(t,x)$.
The PDP associated to the generator $\tA_\infty$ has discontinuities which are
not present for the process $x_C^{N,\epsilon}$. They only appear at the limit
$\epsilon\to 0$. Therefore, the convergence does not occur in $D(\R^+;[0,T])$, indeed
creation of discontinuities can not occur in $D(\R^+;[0,T])$.
As shown below, it is possible to show
tightness, then weak convergence of the process $x_C^{N,\epsilon}$ in $L^p([0,T])$.
Unfortunately, this  topology does not imply a uniform bound on the processes and we cannot
get rid of the boundedness assumptions on the reaction rates as was done in the preceding
sections.

\end{Remark}

{\bf{Proof:}} The proof is divided into four steps.

{\it{Step 1: Limit of $\theta^{N,\epsilon}$}}
\par
We begin with the test-function
$$f(x_C,0)=0, \quad f(x_C,1)=1,\qquad \forall x_C\in\R^{M_C}$$

Then, $\tA_{N, \epsilon}f(x_C,\theta)= \lambda_\theta(x_C,0)\1_{\theta=0} -\frac1\epsilon
\tilde\lambda_\theta(x_C,1)\1_{\theta=1}$, and
$$
\begin{array}{l}
\1_{\theta^{N,\epsilon}(t)=1}-\1_{\theta^{N,\epsilon}(0)=1}-\ds \int_0^t \lambda_\theta(x_C^{N,\epsilon}(s),0)\1_{\theta^{N,\epsilon}(s)=0} -\frac1\epsilon\tilde
\lambda_\theta(x_C^{N,\epsilon}(s),1)\1_{\theta^{N,\epsilon}(s)=1}ds
\end{array}
$$
is a martingale. We then deduce, for every $T>0$ and $t\in [0,T]$,
\begin{equation}
\label{e2.5.3}
\frac{\alpha}{\epsilon} \E\left( \int_0^t \1_{\theta^{N,\epsilon}(s)=1}ds \right) \leq M_\lambda T+1
\end{equation} where $M_\lambda$ is the superior bound of $\lambda_\theta(\cdot,0)$. Then, $\theta^{N,\epsilon}\to 0$ in $L^1(\Omega\times [0,T];\{0,1\})$ when $\epsilon\to0$ and $N\to\infty$.

{\it{Step 2: Tightness of $x_C^{N,\epsilon}$}}
\par

For all $t\in [0,T]$, we write
$$x_C^{N,\epsilon}(t)=x_C^{N,\epsilon}(0)+\sum_{r\in\mathcal R_{DC }}\frac{\gamma_r^C}{N}Y_r\bigg(\displaystyle \int_0^t N\tilde\lambda_r(x_C^{N,\epsilon}(s),0)\1_{\theta^{N,\epsilon}(s)=0}+\frac{N}{\epsilon}\tilde\lambda_r(x_C^{N,\epsilon}(s),1)\1_{\theta^{N,\epsilon}(s)=1} \, ds\bigg)$$ where $(Y_r)_{r\in \mathcal R_{DC }}$ are
independent standard Poisson processes. Then
\begin{equation}
\label{e2.5.4}
\begin{array}{ll}
\left|x_C^{N,\epsilon}(t)\right|&\ds \le \left|x_C^{N,\epsilon}(0)\right|+ \sum_{r\in\mR_{DC }}\frac{ |\gamma_r^C|}{N}
Y_r\left( N M \left(T+
\frac{1}{\epsilon}\int_0^T \un_{\theta^{\epsilon,N}(s)=1} ds\right)
\right)
\end{array}
\end{equation}
where $M$ is the superior bound of $\tilde\lambda_r, r\in \mR_{DC }$.
Since $\sup_{t\in [0,T]}\frac1NY_r(Nt)$ is almost surely bounded, using (\ref{e2.5.3}) and (\ref{e2.5.4}), we have for all $K>0$, for  $N$ large enough and  $\epsilon$ small enough:
$$
\P(\sup_{[0,T]} |x_C^{\epsilon, N}(t)| \ge K) \le \varepsilon(K)
$$
where $\varepsilon (K) \to 0$ if $K\to \infty$.

We introduce $BV(0,T)$ the space of functions of bounded variations on $[0,T]$, with its norm, defined for $f\in BV(0,T)$ by $$\|f\|_{BV(0,T)}=\| f\|_{L^1([0,T])}+\sup\{\sum_i | f(t_{i+1}-f(t_i)|, \textrm{ for } (t_i)_i \textrm{  a finite subdivision of } [0,T]\}$$

Since the processes $Y_r$ are non decreasing, we obtain for a subdivision of $[0,T]$, denoted by $0=t_0<t_1<\dots<t_n=T$,
$$
\begin{array}{l}
\ds \sum_{i=0}^{n-1}|x_C^{\epsilon, N}(t_{i+1}) - x_C^{\epsilon, N}(t_i)|\\
\\
\ds \le \sum_{i=0}^{n-1}\sum_{r\in\mR_{DC }}\frac{| \gamma_r^C|}{N} \bigg[
Y_r\left(\int_0^{t_{i+1}} N\tilde\lambda_r(x_C^{\epsilon, N}(s),0)\un_{\theta^{\epsilon,N}(s)=0} +
\frac{N}{\epsilon}\tilde\lambda_r(x_C^{\epsilon, N}(s),1)\un_{\theta^{\epsilon,N}(s)=1} ds
\right)\\
\\ \ds \hspace{3cm}-Y_r\left(\int_0^{t_i} N\tilde\lambda_r(x_C^{\epsilon, N}(s),0)\un_{\theta^{\epsilon,N}(s)=0} +
\frac{N}{\epsilon}\tilde\lambda_r(x_C^{\epsilon, N}(s),1)\un_{\theta^{\epsilon,N}(s)=1} ds\right)\bigg]\\
\\
\ds = \sum_{r\in\mR_{DC }}\frac{| \gamma_r^C|}{N}
Y_r\left(\int_0^{T} N\tilde\lambda_r(x_C^{\epsilon, N}(s),0)\un_{\theta^{\epsilon,N}(s)=0} +
\frac{N}{\epsilon}\tilde\lambda_r(x_C{\epsilon, N}(s),1)\un_{\theta^{\epsilon,N}(s)=1} ds
\right)
\end{array}
$$

It is then easy to prove that for all $K>0$, for $N$ large enough and  $\epsilon$ small enough:
$$
\P(\|x_C^{N,\epsilon}\|_{BV(0,T)} \ge K) \le \varepsilon(K)
$$
where $\varepsilon (K) \to 0$ if $K\to \infty$.

The set $\{f\in BV(0,T),\textrm{ such that } \| f\|_{BV(0,T)}\leq K\}$ is relatively compact in $L^1([0,T])$ (see for instance \cite{giusti}), and the set $$\{f\in BV(0,T), \textrm{ such that } \| f\|_\infty \leq K \textrm{ and }\| f\|_{BV(0,T)}\leq K\}$$ is relatively compact in $L^p([0,T])$ for $1<p<\infty$. We then conclude that the family of processes $(x_C^{N,\epsilon})_{N,\epsilon}$ is tight in $L^p([0,T])$, $1\leq p<\infty. $

{\it{Step 3: Identification of the limit distribution of $x_C^{N,\epsilon}$}}
\par
Since we are interested by the limit distribution of $x_C^{N,\epsilon}$, we introduce test-functions depending only on the continuous variable $x_C\in\R^{M_C}$. We then define
$$f(x_C,0)=\varphi(x_C),\quad x_C\in\R^{M_C} $$ where $\varphi \in C_b^1(\R^{M_C}). $

We want to define $f(\cdot,1)$ such that
\begin{equation}
\label{e2.5.1.0}
\begin{array}{l}
\ds \sum_{r\in {\mR_{DC}}} \left[ f(x_{C}+\frac1N \gamma_r^C,1)-f(x_{C},1)
\right]N\tilde\lambda_r(x_{C},1)- \tilde\lambda_\theta(x_{C},1) f(x_C,1) \\
\\
\ds =-\tilde\lambda_\theta(x_{C},1)\varphi(x_C).
\end{array}
\end{equation}
Drawing inspiration from the preceding section, we introduce the process $y_N(\cdot,x)$ starting from $x$, and whose generator is
$$
A_N\psi(y)= \sum_{r\in {\mR_{DC}}} \left[ \psi(y+\frac1N \gamma_r^C)-\psi(y)
\right]N\tilde\lambda_r(y,1).
$$
Since $F_1$ is Lipschitz on $\R^{M_C}$, a result of \cite{kurtz71} states that for all $x\in\R^{M_C}$, $T>0$ and $\delta >0$,
\begin{equation}
\label{eqKurtz}
\P\left(\sup_{s\leq T}| y_N(s,x)-\phi_1(s,x)|>\delta\right)\xrightarrow[N \rightarrow \infty]{}0.\end{equation} Remember that $\phi_1(\cdot,x)$ is the flow associated to the vector field $F_1(x)= \sum_{r\in {\mR_{DC}}} \gamma_r^C\tilde\lambda_r(x,1)$ and starting at $x$.

We also introduce the semigroup $\left(P_t^N\right)_{t\geq 0}$ defined on $\mathcal B_b(\R^{M_C})$ by
$$P_t^N\psi(x):=\E\left(\psi\left(y_N(t,x)\right)e^{-\int_0^t\tilde\lambda_\theta\left(y_N(s,x),1\right)ds}\right).$$
It is classical that $(P_t^N)_{t\geq0}$ satisfies the semigroup property and  that
$$\frac{d}{dt}P_t^N\psi(x)=A_NP_t^N\psi(x)-\tilde\lambda_\theta(x,1)P_t^N\psi(x),\quad x\in \R^{M_C},\quad\psi\in \mathcal C_b(\R^{M_C}), t\in\R^+. $$
Then we propose:
$$f(x_C,1): =\int_0^\infty P_s^N\left(\tilde \lambda_\theta(\cdot,1)\varphi\right)\left(x_C \right)ds.$$
Since $\tilde \lambda_\theta(\cdot,1)$ is bounded below by
$\alpha >0$, $f(\cdot,1)$ is well defined. Moreover $f(\cdot, 1)$ satisfies (\ref{e2.5.1.0}). Indeed,
$$
\begin{array}{l}
A_Nf(x_C,1)-\tilde\lambda_\theta(x_C,1)f(x_C,1) \\
={\displaystyle \int_0^\infty A_NP_t^N\left(\tilde\lambda_\theta(\cdot,1)\varphi\right)(x_C)-\tilde\lambda_\theta(x_C,1)P_t^N\left(\tilde\lambda_\theta(\cdot,1)\varphi\right)(x_C)\,dt}\\
={\displaystyle \int_0^\infty \frac{d}{dt}P_t^N\left(\tilde\lambda_\theta(\cdot,1)\varphi\right)(x_C)\,dt}
\\
=\lim_{t\to\infty}P_t^N\left(\tilde\lambda_\theta(\cdot,1)\varphi\right)(x_C)-\left(\tilde\lambda_\theta(\cdot,1)\varphi\right)(x_C)\\
=-\left(\tilde\lambda_\theta(\cdot,1)\varphi\right)(x_C)\\
\end{array}
$$
This test-function $f$ satisfies:
$$
\tA_{N,\epsilon} f(x_C,1) =0
$$
and
$$
\begin{array}{l}
\ds  \tA_{N,\epsilon} f(x_C,0) = \sum_{r\in {\mR_{DC}}} \left[ f(x_{C}+\frac1N \gamma_r^C,0)-f(x_{C},0)
\right]N\tilde\lambda_r(x_{C},0)\\
\\
+\ds   \left[ \int_0^\infty P_t^N \left(\tilde\lambda_\theta(\cdot,1)\varphi\right)(x_C) dt-\varphi(x_{C})
\right]\lambda_\theta(x_{C},0).
\end{array}
$$
Since $f\in C_b(E)$, $f$ is in the domain of $\tA_{N,\epsilon}$, for all $N,\epsilon$.
Then, for all $0\le t_1,\dots,t_\ell \le t$ and
$\psi\in C_b(\R^{M_C\times\ell})$
\begin{equation}
\label{EqMartingale}
\begin{array}{l}
\E\bigg(\bigg[f(x_C^{N,\epsilon}(t),\theta^{N,\epsilon}(t))- f(x_C^{N,\epsilon}(0),\theta^{N,\epsilon}(0)) \\
\\
\ds -\int_0^t  \tA_{N,\epsilon} f(x_C^{N,\epsilon}(s),0)\un_{\theta^{N,\epsilon}(s)=0} ds
\ds \bigg]\psi(x_C^{N,\epsilon}(t_1),\dots,x_C^{N,\epsilon}(t_\ell))\bigg)=0.
\end{array}
\end{equation}

We now are searching a possible formulation for the limit generator. For all $ x\in \R^{M_C}$ and uniformly in $t\in[0,T]$, for all $T>0$, we have
$$
P^N_t \left(\tilde\lambda_\theta(\cdot,1)\varphi\right)(x)\xrightarrow[N \rightarrow \infty]{} \left(\tilde\lambda_\theta(\cdot,1)\varphi\right)(\phi_1(t,x))e^{-\int_0^t \tilde\lambda_\theta (\phi_1(s,x),1)ds}.$$
Indeed,
$$\begin{array}{ll}
\ds P^N_t \left(\tilde\lambda_\theta(\cdot,1)\varphi\right)(x)-\left(\tilde\lambda_\theta(\cdot,1)\varphi\right)(\phi_1(t,x))e^{-\int_0^t \tilde\lambda_\theta (\phi_1(s,x),1)ds}&\\
\ds =\E\left(\left(\tilde\lambda_\theta(\cdot,1)
\varphi\right)(y_N(t,x))e^{-
\int_0^t\tilde\lambda_\theta(y_N(s,x),1)ds}-\left(\tilde\lambda_\theta(\cdot,1)\varphi\right)(\phi_1(t,x))e^{-\int_0^t\tilde\lambda_\theta(y_N(s,x),1)ds}\right)&\\
\ds +\E\left(\left(\tilde\lambda_\theta(\cdot,1)\varphi\right)(\phi_1(t,x))e^{-\int_0^t\tilde\lambda_\theta(y_N(s,x),1)ds}-\left(\tilde\lambda_\theta(\cdot,1)\varphi\right)(\phi_1(t,x))e^{-\int_0^t\tilde\lambda_\theta(\phi_1(s,x),1)ds}\right)\\
\ds=a+b&\\
\end{array}$$
and it can be shown that for all $\delta >0$, and $t\leq T$,
$$\begin{array}{ll}
|a|&\leq \delta\, L_{\tilde\lambda_\theta(\cdot,1)\varphi} +2\,\|\tilde\lambda_\theta(\cdot,1)\varphi\|_\infty\,\P\left(\sup_{s\leq T}| y_N(s,x)-\phi_1(s,x)|>\delta\right)\\
\end{array}$$
and
$$\begin{array}{ll}
|b|&\leq \|\tilde\lambda_\theta(\cdot,1)\varphi\|_\infty\left(T\,\delta\, L_{\tilde\lambda_\theta(\cdot,1)} +2\,\P\left(\sup_{s\leq T}| y_N(s,x)-\phi_1(s,x)|>\delta\right)\right)\\
\end{array}$$
where $L_{\tilde\lambda_\theta(\cdot,1)\varphi}$ and $L_{\tilde\lambda_\theta(\cdot,1)}$ are the Lipschitz constants of the functions $\tilde\lambda_\theta(\cdot,1)\varphi$ and $\tilde\lambda_\theta(\cdot,1)$.
\par
\medskip

By (\ref{eqKurtz})  and  the dominated convergence theorem, we deduce that
$$
\tA_{N} f(x_C,0)\xrightarrow[N \rightarrow \infty]{} \tA_\infty\varphi(x_C).
$$

The tightness in $L^p([0,T])\times L^1(\Omega\times [0,T];\{0,1\})$ of the family of processes $(x_C^{N,\epsilon},\theta^{N,\epsilon})_{N,\epsilon}$, and the Skorohod representation theorem imply the existence of a subsequence which converges almost surely in $L^p([0,T])\times L^1(\Omega\times [0,T];\{0,1\})$ to $(x_C(t),0)_t\in L^p([0,T])\times L^1(\Omega\times [0,T];\{0,1\}) $. This almost surely convergence implies that a new subsequence can be extracted such that
$$
(x_C^{N_k,\epsilon_k}(t),\theta^{N_k,\epsilon_k})  \xrightarrow[k\to\infty]{}( x_C(t),0)
$$
almost surely, and for almost all $t\in [0,T]$.
To simplify notations, this subsequence is still noted $x_C^{N,\epsilon}$.

Taking $N\to\infty$ and $\epsilon\to0$ in (\ref{EqMartingale}), we deduce by the dominated convergence theorem, the boundedness of $f$ and $\psi$, and the hypothesis on $x_C^{N,\epsilon}(0)$ and $\theta^{N,\epsilon}(0)$:
$$
\begin{array}{l}
\E\bigg(\bigg[f(x_C^{N,\epsilon}(t),\theta^{N,\epsilon}(t))- f(x_C^{N,\epsilon}(0),\theta^{N,\epsilon}(0)) \ds \bigg]\psi(x_C^{N,\epsilon}(t_1),\dots,x_C^{N,\epsilon}(t_\ell))\bigg)\\
\\
  \xrightarrow[N\rightarrow\infty,\epsilon \rightarrow 0]{}\ds \E\bigg(\bigg[\varphi(x_C(t))- \varphi(x_C(0)) \ds \bigg]\psi(x_C(t_1),\dots,x_C(t_\ell))\bigg)\\
\end{array}
$$ for almost all $t\in [0,T]$ and almost all $0\leq t_1,\ldots, t_l\leq t$. We now prove that

\begin{equation}
\label{conv}
\tA_{N,\epsilon}f(x_C^{N,\epsilon}(s),0)\1_{\theta^{N,\epsilon}(s)=0}\xrightarrow[N\rightarrow\infty,\epsilon \rightarrow 0]{}\tA_\infty\varphi(x_C(s))
\end{equation}
almost surely, and for almost all $s\in [0,T]$.

First, for $s\in [0,T]$, $t\in\R^+$, we write
$$\begin{array}{ll}
P_t^N\left(\tilde\lambda_\theta(\cdot,1)\varphi\right)(x_C^{N,\epsilon}(s))-\left(\tilde\lambda_\theta(\cdot,1)\varphi\right)\left(\phi_1(t,x_C(s))\right)e^{-\int_0^t\tilde\lambda_\theta(\phi_1(u,x_C(s)),1)du}&\\
&\\
=\E\bigg[\bigg(\tilde\lambda_\theta(\cdot,1)\varphi\bigg)\left(y_N(t,x_C^{N,\epsilon}(s))\right)e^{-\int_0^t\tilde\lambda_\theta\left(y_N\left(u,x_C^{N,\epsilon}(s)\right),1\right)du}\bigg]&\\&\\
-\E\bigg[\bigg(\tilde\lambda_\theta(\cdot,1)\varphi\bigg)\big(\phi_1(t,x_C(s))\big)e^{-\int_0^t\tilde\lambda_\theta\left(y_N\left(u,x_C^{N,\epsilon}(s)\right),1\right)du}\bigg]&\\
&\\
+\E\bigg[\bigg(\tilde\lambda_\theta(\cdot,1)\varphi\bigg)\big(\phi_1(t,x_C(s))\big)e^{-\int_0^t\tilde\lambda_\theta\left(y_N\left(u,x_C^{N,\epsilon}(s)\right),1\right)du}\bigg]&\\
&\\
-\left(\tilde\lambda_\theta(\cdot,1)\varphi\right)\left(\phi_1(t,x_C(s))\right)e^{-\int_0^t\tilde\lambda_\theta(\phi_1(u,x_C(s)),1)du}&\\
&\\
=a_1+b_1&\\
\end{array}$$
For all $\delta>0$, $T_1>0$, and $0\leq t\leq T_1$,
$$\begin{array}{ll}
|a_1|&\leq \delta\, L_{\tilde\lambda_\theta(\cdot,1)\varphi} +2\,\|\tilde\lambda_\theta(\cdot,1)\varphi\|_\infty\,\P\left(\sup_{u\leq T_1}| y_N(u,x_C^{N,\epsilon}(s))-\phi_1(u,x_C(s))|>\delta\right),\\
\end{array}$$
and if $t\geq T_1$,
$$|a_1|\leq 2\|\tilde\lambda_\theta(\cdot,1)\varphi\|_\infty e^{-\alpha T_1}.$$
Similarly, for $0\leq t\leq T_1$,
$$\begin{array}{ll}
|b_1|&\leq  \|\tilde\lambda_\theta(\cdot,1)\varphi\|_\infty\left(T_1\,\delta\, L_{\tilde\lambda_\theta(\cdot,1)} +2\,\P\left(\sup_{u\leq T_1}| y_N(s,x_C^{N,\epsilon}(s))-\phi_1(s,x_C(s))|>\delta\right)\right)\\
\end{array}$$
and if $t\geq T_1$,
$$|b_1|\leq 2\|\tilde\lambda_\theta(\cdot,1)\varphi\|_\infty e^{-\alpha T_1}.$$

We choose $T_1$ sufficiently large and then $\delta$ small enough and have almost surely, for all $s\in [0,T]$, $t\in \R^+$:
$$P_t^N\left(\tilde\lambda_\theta(\cdot,1)\varphi\right)(x_C^{N,\epsilon}(s))-\left(\tilde\lambda_\theta(\cdot,1)\varphi\right)\left(\phi_1(t,x_C(s))\right)e^{-\int_0^t\tilde\lambda_\theta(\phi_1(u,x_C(s)),1)du} \xrightarrow[N\to\infty,\epsilon\to0]{}0.$$
By the dominated convergence theorem, we have almost surely, for almost all $s\in[0,T]$,
$$\displaystyle \int_0^\infty P_t^N\left(\tilde\lambda_\theta(\cdot,1)\varphi\right)(x_C^{N,\epsilon}(s))\, dt \xrightarrow[N\to\infty,\epsilon\to0]{}\displaystyle \int_0^\infty\left(\tilde\lambda_\theta(\cdot,1)\varphi\right)\left(\phi_1(t,x_C(s))\right)e^{-\int_0^t\tilde\lambda_\theta(\phi_1(u,x_C(s)),1)du}\, dt.$$
Since $\varphi\in C^1_b(\R^{M_C})$, we obtain (\ref{conv}).

By the dominated convergence theorem, we conclude finally that
$$
\begin{array}{l}
\E\bigg(\bigg[\varphi(x_C(t))- \varphi(x_C(0))
\ds -\int_0^t  \tA_{\infty} \varphi(x_C(s)) ds
\ds \bigg]\psi(x_C(t_1),\dots,x_C(t_\ell))\bigg)=0
\end{array}
$$ for almost all $t, t_1,\ldots, t_l$ such that $t \in [0,T]$ and $0\leq t_1\ldots t_l\leq t$.

{\it{Step 4: Uniqueness of the solution of the martingale problem}}

By Theorem \ref{t2.2}, proved in a {weaker} sense in the Appendix, we conclude that the limit $x_C$, more exactly a version of $x_C$, is the PDP whose generator is given by $\tilde\mA_\infty$.

%\begin{Remark}
%This model can be extended to the case where the discrete variable $\theta$ takes a finite number of values $\{0, i=1,\ldots,K\}$ such that, if $\theta=0$, the rates of the $\mathcal R_{DC} \cup \mathcal R_{C}$ type reactions \corrO{are of order $N$ or vanish}, and if $\theta\neq 0$, \corrO{some of these rates} are of order $\frac{N}{\epsilon}$. Moreover, the rates of the reactions changing $\theta$ from non-null state to state 0 are of order $\frac1\epsilon$, whereas the rates of the opposite reactions are of order 1.
%\end{Remark}
%JE DOIS ENCORE LE VERIFIER...

\appendix
\section{Appendix}
{\it{Proof of theorem \ref{t2.2}}}
\par

Let $x=(x_t)_{t\geq0}$ denote the canonical process
on $D(\R^+;E)$. Recall that we say that a probability measure $P$ on $D(\R^+;E)$ is solution of the martingale problem associated to the generator $\mathcal A$ if, for all $f\in\mE$,
$$f(x_t)-f(x_0)-\int_0^t\mA f(x_s)ds$$ is a $P$-martingale.

Let $(P_t)_{t\ge 0}$, $P_{x}$ and $\bar \mA$ denote the semigroup, the probability measure and the infinitesimal generator associated to the PDP whose characteristics are $F$, $\lambda$, $Q$, and starting from $x\in E$.

We need the following
\begin{Proposition}
\label{p1.1}
Assume that Hypothesis \ref{h2} holds, then, for all $t\in\R^+$, $f\in \mE$, $P_t f $ is also in $\mE$
and
\begin{equation}
\label{ePt1}
\|P_{t}f\|_{\infty}\le \|f\|_{\infty},
\end{equation}
\begin{equation}
|P_tf(y_1,\nu)-P_tf(y_2,\nu)|\le c  \| f\|_{\mE}e^{Kt} |y_1-y_2|,\quad y_1,\, y_2\in \R^n,\;
\nu \in \N^d,
\end{equation}
for some constants $c$ and $K$ depending only on the characteristics of the PDP.
\end{Proposition}

{\bf Proof :}
Equation (\ref{ePt1}) is clear. Let us define, for $f\in\mE$ and $\psi\in\mB_b(\R^+\times E)$ :

$$
\begin{array}{ll}
G_f\psi(t,x) & \ds= \E_x\left( f(x_t)\un_{t<T_1} +\psi(t-T_1,x_{T_1})\1_{t\ge T_1}\right) \\
\\
& \ds= f(\phi_\nu(t,y),\nu)H(t,x) \\
\\
&\ds+ \int_0^t \int_E \psi(t-s,z)Q(dz;\phi_\nu(s,y),\nu) \lambda(\phi_\nu(s,y),\nu)H(s,x)
ds
\end{array}
$$
for $(t,x)\in \R^+\times E$, with $x=(y,\nu)$. Then, according to Lemma 27.3 of \cite{Davis},
$$
G_f^n\psi(t,x)  \ds= \E_x\left( f(x_t)\1_{t<T_n} +\psi(t-T_n,x_{T_n})\1_{t\ge T_n}\right)
$$
and
$$
\lim_{n\to\infty}G_f^n\psi(t,x)=P_tf(x).
$$
We then deduce by dominated convergence and the strong Markov property:
\begin{eqnarray}
\label{ePt2}
P_tf(x)&= &f(\phi_\nu(t,y),\nu)H(t,x)\nonumber\\
&&+ \int_0^t \int_E P_{t-s}f(z)Q(dz;\phi_\nu(s,y),\nu) \lambda(\phi_\nu(s,y),\nu)H(s,x)ds.
\end{eqnarray}
That is, $P_\cdot f : t\mapsto P_tf$ is a fixed point of $G_f$.

For $T>0$, we introduce the Banach space $\mE_T=L^\infty([0,T],\mE)$, with the norm
$$\|\psi\|_{\mathcal E_T}= \sup_{t\in [0,T]}e^{-\alpha t } \left(\|\psi(t,\cdot)\|_\infty +L_{\psi(t,\cdot)}\right)$$
where $\alpha$ will be fixed hereafter.

For $t\in [0,T]$, $(y,\nu)\in E$, we set $q_t(y,\nu) =f(\phi_\nu(t,y),\nu)H(t,y,\nu)$. Then $q
:\, t\mapsto q_t $ is in $\mathcal E_T$. Indeed, $$
|H(t,y,\nu)|\le 1,\quad (y,\nu)\in E,\; t\in [0,T].
$$
Also, from the Gronwall Lemma, we have that $\phi_\nu$ is differentiable with respect to $y$ and
$$
\left|  D_y\phi_\nu(t,y)  \right| \le e^{Lt} ,\quad (y,\nu)\in E,\;t\in [0,T].
$$
Then we have, for $(y,\nu)\in E$, $t\in [0,T]$,
$$
\begin{array}{ll}
\left|D_yH(t,y,\nu)\right| &\le \ds\int_0^t \left|D_y\{\lambda(\phi_\nu(s,y),\nu)\}\right| ds\le L_\lambda \int_0^t \left|D_y\phi_\nu(s,y)\right| ds\le \frac{L_\lambda}{L} \left(e^{Lt}-1\right).
\end{array}
$$
Hence if $\alpha> L$
$$
\|q\|_{\mathcal E_T}\le \left( L_f+\frac{L_\lambda}L\|f\|_\infty\right)+\|f\|_\infty.
$$
If $\psi\in \mathcal E_T$, we then easily show that $G_f\psi\in \mathcal E_T$. Moreover, if $\psi_1,\; \psi_2 \in \mathcal E_T$, we have
\begin{equation}
\label{e1.2.1}
\| G_f\psi_1(t,\cdot)-G_f\psi_2(t,\cdot)\|_{\infty} \le M_\lambda \int_0^t \|\psi_1(t-s,\cdot)-\psi_2(t-s,\cdot)\|_\infty ds.
\end{equation}
Also, it is easy to see that
$$
\|D_yG_f\psi_1(t,\cdot)-D_yG_f\psi_2(t,\cdot)\|_\infty\\
\le \kappa \int_0^te^{Ls} \|\psi_1(t-s,\cdot)-\psi_2(t-s,\cdot)\|_{\mE}ds
$$
where $\kappa$ is a constant depending only on the characteristics of the PDP.

Then we easily deduce that
$$
\begin{array}{ll}
\|G_f \psi_1-G_f\psi_2\|_{\mathcal E_T} &\ds \le \kappa_1 \sup_{t\in [0,T]} e^{-\alpha t} \left(\int_0^t e^{\alpha(t-s)}ds
+ \int_0^t e^{\alpha (t-s)}e^{Ls}ds\right)\|\psi_1-\psi_2\|_{\mathcal E_T}\\
&\ds\le \kappa_1 \left( \frac1\alpha+\frac1{\alpha-L}\right)\|\psi_1-\psi_2\|_{\mathcal E_T}
\end{array}
$$
where $\kappa_1$ is a constant depending only on the characteristics of the PDP.
Wie now choose $\alpha$ sufficiently large and deduce from Picard theorem that $G_f$ has a unique fixed point in $\mE_T$. This fixed point is the limit $G^n_f\psi$ for any $\psi\in\mE_T$. Thus $P_\cdot f : t\mapsto P_tf$ is that fixed point, and is then in $\mE_T$.

The Lipschitz constant $L_{P_tf}$ of $P_tf$ can easily be found with the preceding results and Gronwall Lemma.

$\square$

The infinitesimal generator  $\bar \mA$ of $P_t$ is characterized by (cf. for instance \cite{ethierkurtz}, part 1, chapter 2):
\begin{equation}
\label{eqInfGene}
(\mu -\bar\mA)^{-1}f=\int_{0}^\infty e^{-\mu t} P_{t}f dt, \quad \mu>0,\quad f\in C_b(E).
\end{equation}

%Marche aussi si on met seulement born√à (JE METS CONTINUE ET BORNEE AU LIEU DE BORNEE, CAR J'UTILISE LA PROPOSITION P. 10 ETHIER ET KURTZ. )
Let $f\in \mE$, $\mu>K$ and $\varphi = (\mu -\bar\mA)^{-1}f$, then $\varphi\in \mathcal D(\bar\mA)$ and
$$
\frac{P_t\varphi(x) -\varphi(x)}{t}  \xrightarrow[t \rightarrow 0]{}  \bar{\mA}\varphi(x),\quad x\in E.
$$

By Proposition \ref{p1.1}, it is easy to see that if Hypothesis \ref{h2} is satisfied, $(\mu -\bar\mA)^{-1}$ maps $\mE$ into itself if $\mu > K$. In particular, $\varphi\in \mE$ and
for all $x\in E$:
\begin{equation}
\label{martingale}
\varphi(x_t)-\varphi(x) -\int_0^t \mA \varphi(x_s)ds
\end{equation}
is a $P_{x}$-martingale.

It follows that, for all $x\in E$, and $t\in \R^+$,
$$
\frac{P_t\varphi(x)-\varphi(x)}{t} =\frac1t \int_0^t\E_x( \mA \varphi(x_s))ds.
$$

Since $\mA \varphi(x_s)$ is $P_x$-a.s. bounded and continuous on $[0,T_1[$, we deduce that

$$
\frac{P_t\varphi(x)-\varphi(x)}t  \xrightarrow[t \rightarrow 0]{}  \mA \varphi(x), \quad x\in E.
$$

It follows that $$\mA\varphi(x)= \bar \mA \varphi(x), \quad x\in E.$$

\begin{Remark}
In fact, we do not need that  \eqref{martingale} defines a $P_x$ martingale. This can be
weakened and replaced by the following:

For all $x\in E$, $\varphi\in \mathcal E$ and for $t\in \mathbb T$, where $\mathbb T$ is a dense subset of $\R^+$:
\begin{equation}
\label{martingale2}
P_t\varphi(x)-\varphi(x) =\int_0^t\E_x\left( \mA \varphi(x_s)\right)ds
\end{equation}
In that case, we can define a sequence $(t_n)_n$ such that $t_n\in\mathbb T$ for all $n$ and $t_n\to 0$ as $n\to \infty$. Taking $n\to\infty $ in the equality
$$
\frac{P_{t_n}\varphi(x)-\varphi(x)}{t_n} =\frac{1}{t_n}\int_0^{t_n}\E_x\left( \mA \varphi(x_s)\right)ds,
$$ we obtain, since $\mA \varphi(x_s)$ is $P_x$-a.s. bounded and continuous on $[0,T_1[$, $$\bar \mA\varphi(x)= \mA \varphi(x), \quad x\in E.$$

\end{Remark}

We now use a classical argument (see for instance
\cite{ethierkurtz}, part 4, chapter 4).
Let $\widetilde{P}_x$ be another solution of the martingale problem. By the same reasoning as before, let $f\in \mE$, $\mu>K$ and $\varphi = (\mu -\bar\mA)^{-1}f$, then
$$
\varphi(x_t)-\varphi(x) -\int_0^t \mA \varphi(x_s)ds
$$
is a $\widetilde{P}_{x}$-martingale.

Let us write
$$
\widetilde {\E}_{x} (\varphi(x_t)-\int_0^t \mA \varphi(x_s)ds ) =\varphi(x).
$$
Multiply this identity by $\mu e^{-\mu t}$ and integrate on $[0,\infty)$ to obtain
%$$
%\widetilde{\E}_{x}\left( \mu \int_0^\infty e^{-\mu t} \varphi(x_t) - \mu \int_0^\infty e^{-\mu t}
%\int_0^t \mA \varphi(x_s)ds \; dt \right)= \varphi(x)
%$$
%and
$$
\widetilde{\E}_{x}\left( \int_0^\infty e^{-\mu t}\left( \mu  \varphi(x_t)- \mA \varphi(x_t)\right) dt
\right) =  \varphi(x).
$$
Since  $\bar \mA\varphi=\mA \varphi$, we deduce, from (\ref{eqInfGene})
$$
\widetilde{\E}_{x}\left( \int_0^\infty e^{-\mu t} f(x_t)dt\right) = \varphi(x)=\int_{0}^\infty e^{-\mu t} P_{t}f(x) dt, \quad \mu>K,
$$
By injectivity of the Laplace transform, this implies
$$
\widetilde {\E}_{x}\left( f(x_t)\right) = P_tf(x), \textrm{for almost all } t\geq 0, \textrm{ for all } x\in E.
$$
We have proved that the law of a solution of the martingale problem at a fixed time $t$ in a dense set of $\R^+$, is uniquely determined. This implies uniqueness for the martingale problem (see \cite{bill} section 14).

\begin{Remark}
Again, the hypothesis that $\tilde \P_x$ is a solution of the martingale problem can be weakened and replaced by the following:
$$\tilde E_x\left(\varphi(x_t)-\varphi(x)-\int_0^t\mA\varphi\left(x_s\right)ds\right) =0, $$
for $ x\in E$ and almost all $t\in\R^+$. The previous calculus are similar, since the Lebesgue integral is null over a negligeable set.
\end{Remark}

%OU SINON : ON PEUT DIRE QUE NOTRE THEO 1.4 UN COROLLAIRE DIRECT DU THEO 4.1  PARTIE 4 DE ETHIER ET KURTZ.

%\bibliography{../bibs/bib_ovidiu,../bibs/markov-bib,../bibs/modules,../bibs/complex,../bibs/Bioinfo,../bibs/gul2,../bibs/radulescu}

\bibliographystyle{apalike}

\end{document}